\documentclass[12pt,oneside,reqno]{amsart}
\usepackage{amssymb}
\usepackage{txfonts}
\usepackage{bbm}
\usepackage{cases}
\usepackage{amsmath}
\usepackage{graphicx}
\usepackage{mathrsfs}
\usepackage{stmaryrd}
\usepackage{amsfonts}
\usepackage{enumerate,amsmath,amssymb,amsthm}

\numberwithin{equation}{section}

\newcommand{\be}{\begin{eqnarray}}
\newcommand{\ee}{\end{eqnarray}}
\newcommand{\ce}{\begin{eqnarray*}}
\newcommand{\de}{\end{eqnarray*}}
\newtheorem{theorem}{Theorem}[section]
\newtheorem{lemma}[theorem]{Lemma}
\newtheorem{remark}[theorem]{Remark}
\newtheorem{definition}[theorem]{Definition}
\newtheorem{proposition}[theorem]{Proposition}
\newtheorem{Examples}[theorem]{Example}
\newtheorem{corollary}[theorem]{Corollary}

\def\v{{\mathrm{v}}}
\def\e{{\mathrm{e}}}
\def\eps{\varepsilon}

\def\p{\partial}

\def\[{{\Big[}}
\def\]{{\Big]}}
\def\<{{\langle}}
\def\>{{\rangle}}
\def\({{\Big(}}
\def\){{\Big)}}

\def\bx{{\mathbf{x}}}

\def\dif{{\mathord{{\rm d}}}}

\def\min{{\mathord{{\rm min}}}}

\def\no{\nonumber}
\def\={&\!\!=\!\!&}
\def\bt{\begin{theorem}}
\def\et{\end{theorem}}
\def\bl{\begin{lemma}}
\def\el{\end{lemma}}
\def\br{\begin{remark}}
\def\er{\end{remark}}

\def\bd{\begin{definition}}
\def\ed{\end{definition}}
\def\bp{\begin{proposition}}
\def\ep{\end{proposition}}
\def\bc{\begin{corollary}}
\def\ec{\end{corollary}}
\def\bx{\begin{Examples}}
\def\ex{\end{Examples}}

\def\cB{{\mathcal B}}

\def\cI{{\mathcal I}}

\def\cL{{\mathcal L}}
\def\cM{{\mathcal M}}

\def\cS{{\mathcal S}}
\def\cT{{\mathcal T}}

\def\mC{{\mathbb C}}

\def\mE{{\mathbb E}}

\def\mI{{\mathbb I}}

\def\mL{{\mathbb L}}
\def\mM{{\mathbb M}}
\def\mN{{\mathbb N}}

\def\mQ{{\mathbb Q}}
\def\mR{{\mathbb R}}
\def\mS{{\mathbb S}}

\def\mZ{{\mathbb Z}}

\def\sB{{\mathscr B}}

\def\sI{{\mathscr I}}
\def\sJ{{\mathscr J}}
\def\sK{{\mathscr K}}
\def\sL{{\mathscr L}}

\def\sP{{\mathscr P}}
\def\sQ{{\mathscr Q}}

\def\geq{\geqslant}
\def\leq{\leqslant}

\allowdisplaybreaks

\begin{document}

\title{$L^p$-maximal hypoelliptic regularity of nonlocal kinetic Fokker-Planck operators}

\date{}

\author{{Zhen-Qing Chen} \ \
and \ \ {Xicheng Zhang}}

\address{Zhen-Qing Chen:
Department of Mathematics, University of Washington, Seattle, WA 98195, USA\\
Email: zqchen@uw.edu
 }
\address{Xicheng Zhang:
School of Mathematics and Statistics, Wuhan University,
Wuhan, Hubei 430072, P.R.China\\
Email: XichengZhang@gmail.com
 }

\begin{abstract}
For $p\in(1,\infty)$, let $u(t,x,\v)$ and  $f(t,x,\v)$ be in $L^p(\mR \times \mR^d \times \mR^d)$
and satisfy the following nonlocal kinetic Fokker-Plank equation on $\mR^{1+2d}$ in the weak sense:
$$
\p_t u+\v\cdot\nabla_x u=\Delta^{{\alpha}/{2}}_\v u+f,
$$
where $\alpha\in(0,2)$ and $\Delta^{{\alpha}/{2}}_\v$ is the usual fractional Laplacian applied to $\v$-variable.
We show that there is a constant $C=C(p,\alpha,d)>0$ such that for any $f(t, x, \v)\in L^p(\mR \times \mR^d \times \mR^d)=L^p(\mR^{1+2d})$,
$$
\|\Delta_x^{{\alpha}/{(2(1+\alpha))}}u\|_p
+\|\Delta_\v^{{\alpha}/{2}}u\|_p\leq C\|f\|_p,
$$
where $\|\cdot\|_p$ is the usual $L^p$-norm in $L^p(\mR^{1+2d}; \dif z)$.
In fact, in this paper the above inequality is established for a large class of  time-dependent  non-local kinetic Fokker-Plank equations on $\mR^{1+2d}$, with $U_t \v$
and $\sL^{\nu_t}_{\sigma_t}$ in place of $\v\cdot \nabla_x$ and $\Delta^{\alpha/2}_\v$. See Theorem \ref{Main} for details.

\end{abstract}

\maketitle

\section{Introduction}

Consider the following classical heat equation in $\mR^{1+d}=\mR \times \mR^d$:
$$
\p_t u=\Delta u+f,
$$
where $\Delta$ is the Laplacian in $\mR^{d}$. It is by now a classical result that for any $p\in(1,\infty)$, there is a constant $C=C(d,p)>0$
such that for all $f (t, x) \in L^p(\mR \times \mR^d)$,
$$
\|\Delta u\|_{L^p(\mR^{1+d})}\leq C\|f\|_{L^p(\mR^{1+d})},
$$
which is an easy consequence of the classical Mihlin's multiplier theorem (cf. \cite{La-So-Ur}),
and plays a basic role in the $L^p$-theory of second-order parabolic equations (cf. \cite{Kr}).
This type of estimate has been extended to the nonlocal L\'evy operators (a class of pseudo-differential operators with non-smooth symbols)
in \cite{Mi-Pr} and \cite{Zh1}.

\medskip

In this paper, we are concerned with
the following kinetic equation in $\mR^{1+2d}$:
\begin{align}\label{Op}
\p_t u+\v\cdot\nabla_xu=\Delta^{\alpha/2}_\v u+f,\quad \alpha\in(0,2],
\end{align}
where $u(t,x,\v)$ and $f(t,x,\v)$ are Borel measurable functions in $\mR^{1+2d}$,
$(t,x,\v)$ stands for the time, position and velocity variables, and $\Delta^{\alpha/2}_\v=-(-\Delta_\v)^{\alpha/2}$
is the usual fractional Laplacian with respect to the velocity variable. When $\alpha=2$, Kolmogorov in \cite{Ko} first constructed the fundamental solution of
degenerate operator $\p_t+\v\cdot\nabla_x-\Delta_\v$.
Observe that, using It\^o's formula, it is easy to verify that
the infinitesimal generator of the diffusion process
$$
 t\mapsto \left(x_0 -\int_0^t X^{\v_0}_s \dif s, X^{\v_0}_t \right)
$$
where $X^{\v_0}_s$ is a Browian motion on $\mR^d$ starting from $\v_0$ with infinitesimal generator $\Delta$,
is $\Delta_\v -\v\cdot \nabla_x$. Thus for $T>0$,
the solution $u(t, x, \v)$ to \eqref{Op} on $(-\infty, T]\times \mR^d \times \mR^d$
with $\alpha=2$ and $u(T, x, \v)=0$ is given by
$$
u(t, x, \v)=
\mE
 \left[ \int_0^{T-t} f \left(T-t-s, x-\int_0^s X^\v_r \dif r, X^\v_s \right) \dif s \right].
$$
In \cite{Ho}, H\"ormander  established a famous hypoelliptic theorem for general second order partial differential operators.
A more precise {\it global} hypoelliptic regularity estimates are established by Bouchut in \cite{Bo} in 2002:
\begin{equation}\label{e:1.2}
\|\Delta_\v u\|_2+\|\Delta_x^{1/3} u\|_2\leq C\|f\|_2.
\end{equation}
Note that for ``nice" $f(t, x, \v)$ on $\mR\times \mR^d\times \mR^d$,
$$
u(t, x, \v):= -\int_t^\infty f(s, x+\v(s-t), \v)ds
$$
is a solution to $\partial_t u +\v\cdot \nabla_x u =f$.
One can show directly  (see \cite{Bo}) that for any $\alpha >0$,
$$
\| \Delta_x^{\alpha/(2(1+\alpha))} u \|_2 \leq c \, \| \Delta_\v^{\alpha/2} u \|_2^{1/(1+\alpha)}\,
\| f\|_2^{\alpha/(1+\alpha )}.
$$
In particular, taking $\alpha=2$ yields
$$
\| \Delta_x^{1/3} u \|_2 \leq c \, \| \Delta_\v u \|_2^{1/3}\, \| f\|_2^{2/3}.
$$
This explains the mystery of 1/3 appeared in the exponent of $\Delta^{1/3}_x$ in \eqref{e:1.2}.
When $p\not=2$,
through establishing
some weak-type $(1,1)$ estimate, Bramanti, Cupini, Lanconelli and Priola \cite{Br-Cu-La-Pr}
proved the following global regularity estimate
$$
\|\Delta_\v u\|_p\leq C\|f\|_p,\ \ p\in(1,\infty),
$$
which, together with a result of Bouchut in \cite{Bo}, also yields that
$$
\|\Delta_x^{1/3} u\|_p\leq C\|f\|_p,\ \ p\in(1,\infty).
$$
It should be noted that the optimal {\it local} $L^p$-estimates for hypoelliptic differential operators
have been studied by Rothschild and Stein in \cite{Ro-St}, where $\| \Delta^{1/3}_x u \|_2$ term first appeared.

\medskip

On the other hand, for $\alpha\in(0,2)$,  Alexander
\cite{Al} proved the following $L^2$-regularity estimate for \eqref{Op} by using Fourier's transformation,
$$
\|\Delta_x^{\alpha/(2(1+\alpha))}u\|_2+\|\Delta_\v^{\alpha/2}u\|_2\leq C\|f\|_2.
$$
A natural question arises
 whether the above fractional hypoellipticity
estimate still holds for general $p\in(1,\infty)$. Clearly, such type estimates belong to the theory of singular integral operators.
In fact, as
 pointed out in \cite{Al}, the main motivation of studying the
above nonlocal regularity
also comes
 from the investigation of spacially inhomogeneous Boltzmann equations.
Let us explain this point in detail (see also \cite{Vi}).
Denote by $\v$ and $\v_*$ the velocities of two particles
immediately before the collision, and $\v'$ and $\v_*'$ their velocities immediately after the collision.
Physics law says
$$
\v'=\v-\<\v-\v_*,\omega\>\omega, \quad
 \v'_*=\v_*+\<\v-\v_*,\omega\>\omega,\  \omega\in\mS^{d-1},
$$
where $\mS^{d-1}$ is the unit sphere in $\mR^d$.
We have the following relations:
$$
\v+\v_*=\v'+\v'_*, \quad  |\v-\v_*|=|\v'-\v'_*|, \quad  |\v|^2+|\v_*|^2=|\v'|^2+|\v'_*|^2,
$$
(i.e. conservation of velocities and conservation of energies) and
$$
 \<\v',\omega\>=\<\v_*,\omega\>,
 \quad   \<\v'_*,\omega\>=\<\v,\omega\>.
$$
Let $f$ be the density of gases.
The classical Boltzmann equation
says
$$
\p_t f(t,x,\v)+\v\cdot \nabla_x f(t,x,\v)=Q(f,f)(t,x,\v),
$$
 where
$Q(f,g)$ is the collision operator defined by
$$
Q(f,g)(\v):=\int_{\mR^d}\int_{\mS^{d-1}}(f(\v'_*)g(\v')-f(\v_*)g(\v))B(|\v-\v_*|,\omega)\dif\omega\dif \v_*,
$$
where
$$
B(|\v-\v_*|,\omega)=|\v-\v_*|^\gamma b(|\<\v-\v_*,\omega\>|/|\v-\v_*|),
$$
and
$$
\mbox{$b(s)\asymp s^{-1-\alpha}$, $\alpha\in(0,2)$ and $\gamma+\alpha\in(-1,1)$,}
$$
where $\asymp$ means that both sides are comparable up to a constant.
Here and below, we drop ``$(t,x)$'' for simplicity.
By an elementary calculation,
the collision operator has the following Carleman's representation (see Appendix 4.1):
\begin{align}\label{Bol}
\begin{split}
Q(f,g)(\v)=2\int_{\mR^d}\!\!\int_{\{h\cdot w=0\}}&\Big[f(\v-h)g(\v+w)-f(\v-h+w)g(\v)\Big]\\
&\times B(|h-w|,w/|w|)|w|^{1-d}\dif h\dif w.
\end{split}
\end{align}
In particular,
when  $b(s)=s^{-1-\alpha}$,  we can split
$Q$ into two parts
$$
Q(f,g)=Q_1(f,g)+Q_2(f,g),
$$
where $Q_1(f,g)(\v):=g(\v) H_f(\v)$ with
$$
H_f(\v):=2\int_{\mR^d}\!\!\int_{\{h\cdot w=0\}}(f(\v-h)-f(\v-h+w))\frac{|h-w|^{\gamma+1+\alpha}}{|w|^{\alpha+d}}\dif h\dif w,
$$
and
$$
Q_2(f,g)(\v):=\int_{\mR^d}(g(\v+w)-g(\v))\frac{K_f(\v,w)}{|w|^{\alpha+d}}\dif w
$$
with
$$
K_f(\v,w):=2\int_{\{h\cdot w=0\}}f(\v-h)|h-w|^{\gamma+1+\alpha}\dif h.
$$
The linearized Boltzmann equation then takes the following form
that involves non-local operator of fractional Laplacian type:
$$
\p_t g+\v\cdot \nabla_x g=
{\rm p.v.} \int_{\mR^d}
(g(\cdot+w)-g(\cdot))\frac{K_f(\cdot,w)}{|w|^{\alpha+d}}\dif w+g\,H_f.
$$
Note that $K_f$ is a symmetric kernel in $w$, i.e.,
$K_f(\cdot,w)=K_f(\cdot,-w)$, and $g H_f$ is a
zero order term in $g$.

\medskip

The goal of this paper is
to study the following nonlocal kinetic Fokker-Planck equation:
$$
\p_s u+U_s\v\cdot\nabla_x u+\lambda u=\int_{\mR^d}\Big[u(\cdot+\sigma_s w)+u(\cdot+\sigma_s w)-2u(\cdot)\Big]\nu_s(\dif w),
$$
where $\lambda\geq 0$, $\nu_s: \mR_+\to\mL^{sym, (\alpha)}_{non}$ and $U_s, \sigma_s:\mR_+\to\mM^d_{non}$ are measurable functions. Here,
$\mL^{sym, (\alpha)}_{non}$ is the space
of non-degenerate symmetric $\alpha$-stable L\'evy measures and $\mM^d_{non}$ is the space of all nonsingular $d\times d$-matrices.
Under suitable assumptions on $\nu,\sigma$ and $U$, we
will establish in Theorem \ref{Main} of this paper
the following
$L^p$-maximal hypoelliptic regularity:
$$
\big\|\Delta_x^{\alpha/(2(1+\alpha))}u^\lambda\big\|_p+\big\|\Delta_\v^{\alpha/2}u^\lambda\big\|_p\leq C\|f\|_p,\ \ p\in(1,\infty).
$$

\medskip

The rest of the paper is organized as follows.
In Section 2, we give some preliminaries.
In particular, we
 derive some estimates about the density of the processes associated with the nonlocal operators. We also recall Fefferman-Stein's theorem.
In Section 3, we prove our main result Theorem \ref{Main}
for $p\not=2$ by
showing the boundedness of suitably defined operators from $L^\infty$ to $BMO$-spaces.
Some useful facts needed in this paper are collected in Subsections 4.1-4.2 of the Appendix of this paper.
The proof of Theorem \ref{Main} for $p=2$ is given in Subsection 4.3. Its proof is new and more elementary
even for the time-independent case (that is, $U_s$ is independent of $s$) studied in Alexander \cite{Al}.
This elementary proof is based on a
direct Fourier transform.

\medskip

Throughout this paper we use the following convention.
The letter $C$ with or without subscripts will denote an unimportant constant, whose value may change in different places.
Moreover, $f\preceq g$
means that $f\leq C g$ for some constant $C>0$, and $f\asymp g$ means that  $C^{-1}g\leq f \leq C g$ for some $C>1$.
We use $:=$ as a way of definition.
For two real numbers $a$ and $b$, $a\vee b:= \max \{a, b\}$, $a \wedge b:=\min \{a, b\}$ and
$a^+:=\max \{a, 0\}$.

\section{Preliminaries}

Let $\mL^{sym}$ be the set of all symmetric L\'evy measures $\nu$ on $\mR^d$, that is,
(positive) measures $\nu$ on $\mR^d$ such that
$$
\nu(-\dif x)=\nu(\dif x),\ \ \nu(\{0\})=0,\ \ \int_{\mR^d}\big(1\wedge|x|^2\big)\nu(\dif x)<+\infty.
$$
We equip $\mL^{sym}$
with the weak convergence topology.
For $\alpha\in(0,2)$, let $\mL^{sym,(\alpha)}\subset\mL^{sym}$ be the set of all symmetric $\alpha$-stable measures $\nu^{(\alpha)}$ with form
\begin{align}
\nu^{(\alpha)}(A)=\int^\infty_0\left(\int_{\mS^{d-1}}\frac{1_A (r\theta)\Sigma(\dif\theta)}{r^{1+\alpha}}\right)\dif r,\quad A\in\sB(\mR^d),\label{Eq4}
\end{align}
where $\Sigma$ is a finite symmetric measure over the sphere $\mS^{d-1}$ (called spherical measure of $\nu^{(\alpha)}$).

We introduce the following notions.
\bd
\begin{enumerate}[(i)]
\item A symmetric $\alpha$-stable measure $\nu^{(\alpha)}\in \mL^{sym,(\alpha)}$ is called non-degenerate if
\begin{align}
\int_{\mS^{d-1}}|\theta_0\cdot\theta|^\alpha\Sigma(\dif\theta)
>0 \quad \hbox{for every } \theta_0\in\mS^{d-1}.   \label{Spe1}
\end{align}
The set of all non-degenerate symmetric $\alpha$-stable measures is denoted by $\mL^{sym,(\alpha)}_{non}$.
\item For $\nu_1,\nu_2\in\mL^{sym}$, we say that $\nu_1$ is less than $\nu_2$ (simply written as $\nu_1\leq\nu_2$) if
$$
\nu_1( A )\leq \nu_2(A) \quad  \hbox{for any }  A\in\sB(\mR^d).
$$
\end{enumerate}
\ed

\medskip

\br
\rm In this paper, for simplicity we only consider symmetric stable L\'evy measures. This assumption is not crucial.
All the results of this paper can be extended to non-symmetric stable L\'evy measures.
\er

\medskip

For a function $f\in C^2_b(\mR^{d})$, we define the difference operators of first and second orders as follows:
for $x,y\in\mR^d$,
\begin{equation}\label{NB5}
\begin{split}
\delta^{(1)}_{x}f(y):=f(y+x)-f(y),\quad \delta^{(2)}_{x} f(y):=\delta^{(1)}_{x}f(y)+\delta^{(1)}_{-x}f(y).
\end{split}
\end{equation}
Using the fact that
$$ f(y+x)-f(y) = x \cdot \int_0^1 \nabla f (y+sx)\dif s,
$$
we have
for any $p\in[1,\infty]$ and $f\in C^2_b(\mR^d)\cap L^p(\mR^d)$ that
\begin{align}
&\ \|\delta^{(1)}_{x} f\|_p\leq (\|\nabla f\|_p|x|)\wedge(2\|f\|_p),\label{EV311}\\
&\|\delta^{(2)}_{x} f\|_p\leq (2\|\nabla^2 f\|_p|x|^2)\wedge(4\|f\|_p).\label{EV31}
\end{align}

Let $\mM^d$ be the space of all real $d\times d$-matrices and $\mM^d_{non}$ the set of all non-singular matrices.
The identity matrix is denoted by $\mI$, and the transpose of a matrix $\sigma$ is denoted by $\sigma^*$.
Let $\cS(\mR^d)$  be the space of rapidly decreasing functions.
For given $\nu\in\mL^{sym}$, $\sigma\in\mM^d$ and $\alpha\in(0,2)$, we consider the following L\'evy operator:
\begin{align}
\sL^\nu_{\sigma} f(y):=\int_{\mR^d}\delta^{(2)}_{\sigma x} f(y)\nu(\dif x),\ \ f\in\cS(\mR^d).\label{Def1}
\end{align}
Let $\psi^\nu_\sigma$ be the symbol of operator $\sL^\nu_\sigma$, i.e.,
$$
\widehat{\sL^\nu_\sigma f}(\xi)=-\psi^\nu_\sigma(\xi) \hat f(\xi),
$$
where $\hat{f}$ denotes the Fourier transform of $f$. The function $\psi^\nu_\sigma$ is also called a
Fourier multiplier.
It is easy to see that
\begin{align}
\psi^\nu_\sigma(\xi):=2\int_{\mR^d}(1-\cos\<\xi,\sigma x\>)\nu(\dif x).\label{psi}
\end{align}
In particular, for given $\nu_1,\nu_2\in\mL^{sym}$, if $\nu_1\leq\nu_2$, then for any $\sigma\in\mM^d$,
\begin{align}
\psi^{\nu_1}_\sigma(\xi)\leq\psi^{\nu_2}_\sigma(\xi),\ \ \forall\xi\in\mR^d,\label{HG1}
\end{align}
and by \eqref{psi}, \eqref{Eq4} and \eqref{Spe1}, for any $\nu^{(\alpha)}\in\mL^{sym,(\alpha)}_{non}$,
\begin{align}
\psi^{\nu^{(\alpha)}}_\sigma(\xi)\asymp |\sigma^*\xi|^\alpha,\ \xi\in\mR^d,\ \sigma\in\mM^d.\label{HG2}
\end{align}
Moreover, if $\nu(\dif y)=|y|^{-d-\alpha}\dif y$, then $\psi^\nu_\sigma(\xi)=c_{d,\alpha}|\sigma^*\xi|^\alpha,$ where $c_{d,\alpha}$ is a constant only depending on $d,\alpha$.
In this case,
\begin{align}
\sL^\nu_\mI f(y)=c_{d,\alpha}\Delta^{\frac{\alpha}{2}}f(y),\label{Es7}
\end{align}
where $\Delta^{\frac{\alpha}{2}}$ is the usual fractional Laplacian.
In this paper, up to a constant multiple,
we always use the following definition of fractional Laplacian:
\begin{align}
\Delta^{{\alpha}/{2}}f(y)=\int_{\mR^d}\delta^{(2)}_x f(y)\frac{\dif x}{|x|^{d+\alpha}}.\label{Def3}
\end{align}

We have the following commutator estimate.
\bl
Let $\alpha\in(0,2), \sigma\in\mM^d$ and $\nu\in \mL^{sym}, \nu^{(\alpha)}\in\mL^{sym, (\alpha)}$ with
$\nu\leq\nu^{(\alpha)}$. For any $p,q\in[1,\infty]$ with $p\leq q$
 and $\gamma \in ( (\alpha-1)^+,1)$, and for any $\phi\in C^\infty_c(\mR^{d})$, there is a positive constant $C_\phi$ depending on
 $\| \nabla_x^2 \phi \|_p + \| \phi\|_p$
 and $\|\sigma\|, \nu^{(\alpha)},d,\alpha, p,q,\gamma$ such that for any measurable function $f$ on $\mR^{d}$,
\begin{align}
\big\|\sL^{\nu}_{\sigma}(f\phi)-(\sL^{\nu}_{\sigma}f)\phi\big\|_p\leq C_\phi\big([f]_{q,\gamma}+\|f\|_q\big),\label{EV441}
\end{align}
where $[f]_{q,\gamma}:=\sup_x\big(\|\delta^{(1)}_xf\|_q/|x|^\gamma\big)$.
\el
\begin{proof}
By definition \eqref{Def1}, we have
$$
\sL^{\nu}_{\sigma}(f\phi)-(\sL^{\nu}_{\sigma} f)\phi-f\sL^{\nu}_{\sigma}\phi=2\int_{\mR^d}\delta^{(1)}_{\sigma x} f\,\delta^{(1)}_{\sigma x} \phi\nu(\dif x).
$$
Hence, by H\"older's inequality with $\frac{1}{p}=\frac{1}{q}+\frac{1}{r}$ and $\nu\leq\nu^{(\alpha)}$,
\begin{align*}
\|\sL^{\nu}_{\sigma}(f\phi)-(\sL^{\nu}_{\sigma} f)\phi\|_p\leq\|f\|_q\,\|\sL^{\nu}_{\sigma}\phi\|_r
+2\int_{\mR^d}\|\delta^{(1)}_{\sigma x} f\|_q\,\|\delta^{(1)}_{\sigma x} \phi\|_r\nu^{(\alpha)}(\dif x).
\end{align*}
Notice that
$$
\|\delta^{(1)}_{\sigma x} f\|_q\leq\big ([f]_{q,\gamma}|\sigma x|^\gamma\big)\wedge (2\|f\|_q).
$$
The desired estimate then follows by \eqref{Eq4}, \eqref{EV311} and \eqref{EV31}.
\end{proof}

\subsection{Fundamental solutions of nonlocal kinetic Fokker-Planck operator}

In the following, for a function $f(x,\v)\in C^2_b(\mR^{2d})$, we shall write
\begin{align*}
\delta^{(1)}_\v f(x,\v'):=\delta^{(1)}_\v f(x,\cdot)(\v'),& \quad \delta^{(2)}_\v f(x,\v'):=\delta^{(2)}_\v f(x,\cdot)(\v'),\\
\sL^{\nu}_{\sigma,\v}f(x,\v):=\sL^\nu_\sigma f(x,\cdot)(\v),& \quad \Delta^{\frac{\alpha}{2}}_\v f(x,\v):=\Delta^{\frac{\alpha}{2}} f(x,\cdot)(\v),
\end{align*}
and similarly for $\delta^{(1)}_x f(x',\v),\ \delta^{(2)}_x f(x',\v),\ \Delta^{\frac{\alpha}{2}}_x f(x,\v).$

Let $\sigma, U: \mR\to\mM^d_{non}$ be two matrix-valued measurable functions with
\begin{align}
\left\{
\begin{aligned}
&\kappa_0:=\|\sigma\|_\infty+\|\sigma^{-1}\|_\infty+\|U\|_\infty+\sup_{s<t}\Big((t-s)\|\Pi_{s,t}^{-1}\|\Big)<\infty,\\
&\mbox{ where }\ \ \Pi_{s,t}:=\int^t_s U_r\dif r,\ s,t\in\mR \mbox{ with }s<t.
\end{aligned}\label{Kapp}
\right\}
\end{align}
The above assumptions correspond to the non-degeneracy on $\sigma$ and $U$.
Let $\nu:\mR\to\mL^{sym}$ be a measurable map and satisfying that for some $\alpha\in(0,2)$,
\begin{align}
\nu^{(\alpha)}_1\leq\nu_s\leq\nu^{(\alpha)}_2,\quad \nu^{(\alpha)}_1, \nu^{(\alpha)}_2\in\mL^{sym,(\alpha)}_{non}.\label{Con2}
\end{align}
Notice that by \eqref{HG1} and \eqref{HG2}, there is a constant $\kappa_1\in(0,1)$ depending on $\kappa_0$ and $\alpha$ such that
\begin{align}
\kappa_1|\xi|^\alpha\leq \psi^{\nu_s}_{\sigma_s}(\xi)\leq \kappa^{-1}_1|\xi|^\alpha,\ \ \xi\in\mR^d.\label{Con1}
\end{align}
By the above notations, we consider the following time-dependent nonlocal kinetic Fokker-Planck operator
\begin{align}
\sK_s f(x, \v):=\sL^{\nu_s}_{\sigma_s,\v} f(x, \v)+(U_s\v\cdot\nabla_x) f(x,\v).\label{Eq218}
\end{align}

In this subsection we study the existence of smooth fundamental solutions for $\sK_s$ by using
a probabilistic approach, and establish some
short time asymptotic estimates for the heat kernel. Note that the existence of smooth fundamental solution of nonlocal H\"ormander operators was studied
in \cite{Zh2, Zh3, Zh4} (see also the references therein).

Let $N(\dif t,\dif \v)$ be the Poisson random measure on $\mR^{1+d}$
with intensity measure $\nu_t(\dif \v)\dif t$, and $\tilde N(\dif t,\dif\v):=N(\dif t,\dif \v)-\nu_t(\dif \v)\dif t$ the compensated Poisson random measure.
For $s\leq t$, define
\begin{align}
L_{s,t}:=\int^t_s\!\!\!\int_{|\v|\leq 1}\v\tilde N(\dif r,\dif \v)+\int^t_s\!\!\!\int_{|\v|>1}\v N(\dif r,\dif \v),\label{ES1}
\end{align}
and $\Pi_{s,t}:=\int^t_s U_r\dif r$ as well as
\begin{align}
\begin{split}\label{KT0}
K_{s,t}:=(X_{s,t},V_{s,t}):&=\left(\int^t_s U_r\left[\int^r_s\sigma_{r'}\dif L_{s,r'}\right]\dif r,\int^t_s\sigma_r\dif L_{s,r}\right)\\
&=\left(\int^t_s \Pi_{r,t}\sigma_r\dif L_{s,r},\int^t_s\sigma_r\dif L_{s,r}\right),
\end{split}
\end{align}
where the second equality is due to Fubini's theorem.
Notice that $(X_{s,t},V_{s,t})$ solves the following liner SDE:
\begin{align}\label{SDE}
\dif (X_{s,t},V_{s,t})=(U_tV_{s,t},0)\dif t+(0,\sigma_t\dif L_{s,t}),\  (X_{s,s},V_{s,s})=(0,0),\ \ t\geq s.
\end{align}
For any $s\leq t$ and $x,\v\in\mR^d$, let
$$
K_{s,t}(x,\v):=K_{s,t}+(x+\Pi_{s,t}\v,\v)=(X_{s,t}+x+\Pi_{s,t}\v, V_{s,t}+\v),
$$
which solves \eqref{SDE} with starting point $(x,\v)$. In particular,
$$\{K_{s,t}(x,\v), t\geq s, (x,\v)\in\mR^{2d}\}$$ forms a family of time-inhomogenous Markov processes. Let $\cT_{s,t}$ be the associated Markov operator:
\begin{align}
\cT_{s,t} f(x,\v):=\mE f(K_{s,t}(x,\v)),\ \ f\in \cB_b(\mR^d),\label{TST}
\end{align}
where $\cB_b(\mR^d)$ is the set of bounded measurable functions on $\mR^d$.
Clearly, for each $t\geq s$ and $p\in[1,\infty]$,  $\cT_{s,t}$ is a contraction operator in $L^p(\mR^{2d})$ and
\begin{align}
\cT_{s,t}f=\cT_{s,r}\cT_{r, t}f,\ \ s\leq r\leq t.\label{Sem}
\end{align}
Moreover, for any $f\in C^2_b(\mR^{2d})$, $\cT_{s,t}f$ satisfies the following backward Kolmogorov's equation (for example, see \cite{Zh0}):
for Lebesgue-almost all $s\leq t$ and all $x,\v\in\mR^d$,
\begin{align}
\p_s\cT_{s,t} f(x,\v)+\sK_s\cT_{s,t} f(x,\v)=0,\label{EQ00}
\end{align}
where $\sK_s$ is defined by \eqref{Eq218}.
The Fourier transform of $\cT_{s,t}f$ is given by
\begin{align}\label{TST0}
\widehat{\cT_{s,t} f}(\xi,\eta)=\mE\e^{-{\rm{i}}\<(\xi,\eta-\Pi_{s,t}^*\xi), K_{s,t}\>}\hat f(\xi,\eta-\Pi^*_{s,t}\xi).
\end{align}

Below, we use the following convention: If a quantity depends on $\nu$, $\sigma$ and $U$, and when we want to emphasize the dependence,
we shall write them in the place of superscript. For example, there is no further declarations, we sometimes use
$X^{\sigma, U}_{s,t}$, $V^\sigma_{s,t}$,  $K^{\nu}_{s,t}$, $\cT^{\nu}_{s,t}$, and so on.

First of all, we have
\bl
Under \eqref{Con2}, for any $q\in[0,\alpha)$, there is a constant $C=C(d,\nu^{(\alpha)}_2, q,\alpha)>0$
such that for any bounded measurable function $f:\mR\to\mM^d$ and $s<t$,
\begin{align}\label{Mom}
\mE\left|\int^t_s f_r\dif L^\nu_{s,r}\right|^q\leq C\|f\|^q_{L^\infty(s,t)}(t-s)^{\frac{q}{\alpha}}.
\end{align}
\el
\begin{proof}
Since $\nu$ is symmetric, we can write
$$
L_{s,t}=\int^t_s\!\!\!\int_{|\v|\leq(t-s)^{1/\alpha}}\v\tilde N(\dif r,\dif \v)+\int^t_s\!\!\!\int_{|\v|>(t-s)^{1/\alpha}}\v N(\dif r,\dif \v).
$$
By H\"older's inequality and the isometry of stochastic integral, we have
\begin{align*}
\mE\left|\int^t_s\!\!\!\int_{|\v|\leq(t-s)^{1/\alpha}}f_r\v\tilde N(\dif r,\dif \v)\right|^q
&\leq\left(\mE\left|\int^t_s\!\!\!\int_{|\v|\leq(t-s)^{1/\alpha}}f_r\v\tilde N(\dif r,\dif \v)\right|^2\right)^{q/2}\\
&=\left(\int^t_s\!\!\!\int_{|\v|\leq(t-s)^{1/\alpha}}|f_r\v|^2\nu_r(\dif \v)\dif r\right)^{q/2}\\
&\!\!\stackrel{\eqref{Con2}}{\leq} \|f\|_{L^\infty(s,t)}^q\left((t-s)\int_{|\v|\leq(t-s)^{1/\alpha}}|\v|^2\nu^{(\alpha)}_2(\dif \v)\right)^{q/2}\\
&\stackrel{\eqref{Eq4}}{\preceq} \|f\|^q_{L^\infty(s,t)}(t-s)^{\frac{q}{\alpha}}.
\end{align*}
If $q\in(1,\alpha)$, then by Burkholder's inequality (see \cite[(2.10)]{So-Zh}),
\begin{align*}
\mE\left|\int^t_s\!\!\!\int_{|\v|>(t-s)^{1/\alpha}}f_r\v N(\dif r,\dif \v)\right|^q
&\preceq\left(\int^t_s\!\!\!\int_{|\v|>(t-s)^{1/\alpha}}|f_r\v| \nu_r(\dif \v)\dif r\right)^q\\
&\quad+\int^t_s\!\!\!\int_{|\v|>(t-s)^{1/\alpha}}|f_s\v|^q \nu_r(\dif \v)\dif r\\
&\leq\|f\|^q_{L^\infty(s,t)}\left((t-s)\int_{|\v|>(t-s)^{1/\alpha}}|\v| \nu^{(\alpha)}_2(\dif \v)\right)^q\\
&\quad+\|f\|^q_{L^\infty(s,t)}(t-s)\int_{|\v|>(t-s)^{1/\alpha}}|\v|^q \nu^{(\alpha)}_2(\dif \v)\\
&\stackrel{\eqref{Eq4}}{\preceq} \|f\|^q_{L^\infty(s,t)}(t-s)^{\frac{q}{\alpha}}.
\end{align*}
If $q\in(0,1]$, then
\begin{align*}
&\mE\left|\int^t_s\!\!\!\int_{|\v|>(t-s)^{1/\alpha}}f_r\v N(\dif r,\dif \v)\right|^q
\leq \mE\left(\int^t_s\!\!\!\int_{|\v|>(t-s)^{1/\alpha}}|f_s\v|^q N(\dif r,\dif \v)\right)\\
&\qquad=\int^t_s\!\!\!\int_{|\v|>(t-s)^{1/\alpha}}|f_r\v|^q \nu_r(\dif \v)\dif r
\stackrel{\eqref{Eq4}}{\preceq} \|f\|^q_{L^\infty(s,t)}(t-s)^{\frac{q}{\alpha}}.
\end{align*}
Combining the above calculations, we obtain the desired estimate.
\end{proof}

The following is a crucial lemma of this paper.

\bl
Under \eqref{Kapp} and \eqref{Con2},
the random variable $K^{\nu}_{s,t}$ defined by \eqref{KT0} has a smooth density $p^{\nu}_{s,t}(x,\v)$.
Moreover,  for any $n,m\in\mN_0$ and $q_1,q_2\in[0,\alpha)$ with $q_1+q_2<\alpha$,
there exists a positive constant $C=C(d,n,m,\kappa_0,\nu^{(\alpha)}_i,q_i,\alpha)$ such that for all $s<t$,
\begin{align}
\int_{\mR^{2d}}|x|^{q_1}|\v|^{q_2}|\nabla^{n}_{x}\nabla^{m}_{\v} p^{\nu}_{s,t}(x,\v)|\dif x\dif \v\leq
C(t-s)^{((q_1-n)(1+\alpha)+q_2-m)/{\alpha}}.\label{EV11}
\end{align}
\el

\begin{proof}
We divide the proof into four steps. All the constants below will depend only on $d,n,m,\kappa_0,\nu^{(\alpha)}_i,q_i,\alpha$.
\medskip\\
{\bf (i)} First of all, we assume that
$$
\nu_s=\nu^{(\alpha)}\in\mL^{sym,(\alpha)}_{non}.
$$
Let $L^{\nu^{(\alpha)}}$ be an $\alpha$-stable process with the L\'evy measure $\nu^{(\alpha)}$.
Since for any $s<t$, $L^{\nu^{(\alpha)}}$ has the following scaling property:
$$
\Big(L^{\nu^{(\alpha)}}_{(t-s)r}\Big)_{r\geq 0}\stackrel{(d)}{=} \Big((t-s)^{\frac{1}{\alpha}}L^{\nu^{(\alpha)}}_r\Big)_{r\geq 0},
$$
by \eqref{KT0} and the change of variables, we have
\begin{equation}\label{NB4}
\begin{split}
K^{\nu^{(\alpha)}}_{s,t}&\stackrel{(d)}{=}\left((t-s)^{\frac{1}{\alpha}+1}\int^1_0\Pi^{\tilde U}_{r,1}\tilde\sigma_r\dif L^{\nu^{(\alpha)}}_r,
(t-s)^{\frac{1}{\alpha}}\int^1_0\tilde\sigma_{r}\dif L^{\nu^{(\alpha)}}_r\right)\\
&=\Big((t-s)^{\frac{1}{\alpha}+1}X^{\tilde\sigma, \tilde U}_{0,1}, (t-s)^{\frac{1}{\alpha}}V^{\tilde\sigma}_{0,1}\Big),
\end{split}
\end{equation}
where $\tilde U_r:=U_{(t-s)r+s}$ and $\tilde\sigma_r:=\sigma_{(t-s)r+s}$.
This implies that
\begin{align}
p^{\nu^{(\alpha)},\sigma,U}_{s,t}(x,\v)=(t-s)^{-\frac{2d}{\alpha}-d}
p^{\nu^{(\alpha)},\tilde\sigma, \tilde U}_{0,1}((t-s)^{-\frac{1}{\alpha}-1}x,(t-s)^{-\frac{1}{\alpha}}\v).\label{NB3}
\end{align}
Hence, if one can show that for any $n,m\in\mN_0$,
\begin{align}
\int_{\mR^{2d}}|x|^{q_1}|\v|^{q_2}|\nabla^n_x\nabla_\v^mp^{\nu^{(\alpha)},\tilde\sigma,\tilde U}_{0,1}(x,\v)|\dif x\dif\v\leq C,\label{LK3}
\end{align}
then \eqref{EV11} for $\nu_s=\nu^{(\alpha)}$ immediately follows by \eqref{NB3}.
\medskip\\
{\bf (ii)} We make the following further decomposition:
\begin{align}
\nu^{(\alpha)}=\nu_1+\nu_2,\ \ \nu_{1}(\dif \v):=\nu^{(\alpha)}(\dif\v) 1_{|\v|\leq 1},\ \nu_{2}(\dif\v):=\nu^{(\alpha)}(\dif\v) 1_{|\v|>1}.
\end{align}
Let $L^{\nu_i}, i=1,2$ be two independent L\'evy processes with the L\'evy measures $\nu_i$ respectively.
For $i=1,2$, let $K^{\nu_i}_{s,t}=(X^{\nu_i}_{s,t},V^{\nu_i}_{s,t})$ be defined as in \eqref{KT0} with  $\tilde \sigma, \tilde U$
and $L^{\nu_i}_{t-s}$ in place of $\sigma, U$ and $L^\nu_{s,t}$.
In particular,
\begin{align}
K^{\nu^{(\alpha)},\tilde\sigma,\tilde U}_{s,t}\stackrel{(d)}{=}K^{\nu_1}_{s,t}+K^{\nu_2}_{s,t},\label{LK2}
\end{align}
which implies that
$$
p_{0,1}^{\nu^{(\alpha)},\tilde\sigma,\tilde U}(x,\v)=\mE p_{0,1}^{\nu_1}\Big(x-X^{\nu_2}_{0,1},\v-V^{\nu_2}_{0,1}\Big),
$$
where $p_{0,1}^{\nu_1}(x,\v)$ is the distributional density of $K^{\nu_1}_{0,1}$. In view of $q_1+q_2<\alpha$, by \eqref{KT0} and \eqref{Mom}, we have
$$
\mE \left[ (1+|X^{\nu_2}_{0,1}|^{q_1})(1+|V^{\nu_2}_{0,1}|^{q_2})\right]<\infty.
$$
Thus, in order to show \eqref{LK3}, it suffices to prove that for any $n,m\in\mN_0$,
\begin{align}
\int_{\mR^{2d}}(1+|x|^{q_1})(1+|\v|^{q_2})|\nabla^n_x\nabla_\v^mp^{\nu_1}_{0,1}(x,\v)|\dif x\dif\v\leq C.\label{LK33}
\end{align}
\medskip\\
{\bf (iii)} Below, for simplicity of notation, we drop the tilde over $\tilde\sigma, \tilde U$.
Recall for $s<t$,
$
K^{\nu_1}_{s, t}=\left(\int^t_s\Pi_{r, t}\sigma_r \dif  L^{\nu_1}_r, \int^t_s\sigma_r\dif L^{\nu_1}_r\right).
$
By step function approximation, we have
\begin{align}
\mE\e^{\mathrm{i}\<(\xi,\eta), K^{\nu_1}_{s, t}\>}
&=\mE\exp\left( \mathrm{i}\int^t_s\Big<\sigma^*_r\Pi^*_{r, t}\xi+\sigma^*_r\eta, \dif L^{\nu_1}_r\Big>\right) \no\\
 &= \exp
\left( -\int^t_s\psi^{\nu_1}\left(\sigma^*_r \Pi^*_{r, t}\xi+\sigma^*_r\eta\right)\dif r\right),
\label{NB9}
\end{align}
where $\psi^{\nu_1}$ is the characteristic exponent of $L^{\nu_1}$, that is, $\mE \e^{{\rm i}\xi L^{\nu_1}_t}=\e^{-t\psi^{\nu_1}(\xi)}$,
which has the following expression
\begin{align}
\psi^{\nu_1}(\xi)=\int_{|\v|\leq 1}(1-\cos\<\xi, \v\>)\nu^{(\alpha)}(\dif \v).\label{FG1}
\end{align}
Denote the L\'evy exponent of $L^{\nu^{(\alpha)}}$ and $L^{\nu_2}$ by 
$\psi^{\nu^{(\alpha)}}$ 
and $\psi^{\nu_2}$, respectively. Then $\psi^{\nu_2}$ is bounded and 
$\psi^{\nu^{(\alpha)}}(\xi)\asymp |\xi|^\alpha$. 
Hence there are constants $M\geq 1$ and
$c_0>0$ so that
$$ \psi^{\nu_1}(\xi)= 
\psi^{\nu^{(\alpha)}}(\xi)
-\psi^{\nu_2} (\xi) \geq c_0 |\xi|^\alpha
\quad \hbox{for } |\xi| \geq M.
$$
On the other hand, note that
$$
\psi^{\nu_1}(\xi)\geq \int_{|\v|\leq 1/(2M)}(1-\cos\<\xi, \v\>)\nu^{(\alpha)}(\dif \v)=:\psi^{\nu_1}_1 (\xi).
$$
Since $\psi^{\nu_1}_1$ is $C^\infty$-smooth with $\nabla \psi^{\nu_1}_1 (0) =0$,
we have  
\begin{align*}
\psi^{\nu_1}_1 (\xi)&=\int^1_0\!\!\!\int^1_0\!\!\int_{|\v|\leq1/(2M)}\<\xi,\v\>^2\cos\<ss'\xi,\v\>\nu^{(\alpha)}(\dif \v)\dif s\dif s'\\
&\geq\cos(\tfrac{1}{2})\int_{|\v|\leq1/(2M)}\<\xi,\v\>^2\nu^{(\alpha)}(\dif \v)\geq c_0|\xi|^2 \quad \hbox{for } |\xi | \leq M.
\end{align*}
Thus by decreasing the value of $c_0$ if needed, we have
\begin{equation}\label{e:2.34}
\psi^{\nu_1} (\xi) \geq c_0 (|\xi|^2 \wedge |\xi|^\alpha ) \quad \hbox{for } \xi \in \mR^d.
\end{equation}
Hence
\begin{align}
&\int^t_s \psi^{\nu_1} \left(\sigma^*_{r}\Pi^*_{r,t}\xi+\sigma^*_{r} \eta\right)\dif s\geq c_0\int^t_s
\left|\sigma^*_{r} \Pi^*_{r, t}\xi+\sigma^*_{r} \eta\right|^2 \wedge
\left|\sigma^*_{r} \Pi^*_{r, t}\xi+\sigma^*_{r} \eta\right|^\alpha
\dif s\no\\
&\qquad\geq \frac{c_0 (\|\sigma^{-1}\|^{-\alpha}_\infty \wedge \|\sigma^{-1}\|^{-2}_\infty) }
{(\kappa_0+1)^{2-\alpha}} \, \left( |((t-s) \xi,\eta)|^2 \wedge |((t-s) \xi,\eta)|^\alpha \right) \no \\
& \quad \times \inf_{|\bar\xi|^2+|\bar\eta|^2=1}
\int^t_s\left|(t-s)^{-1} \Pi^*_{r, t}\bar\xi+\bar\eta\right|^2 \dif r. \label{NB8}
\end{align}
Fix $\delta,\eps\in(0,1/2)$ being small, whose values will be determined below.
For $\bar\xi,\bar\eta\in\mR^{2d}$ with $|\bar\xi|^2+|\bar\eta|^2=1$, we have either $|\bar\xi|^2\geq 1-\delta$ or $|\bar\eta|^2\geq \delta$.
Since $|a+b|^2 \geq\frac{1}{2}|a|^2-|b|^2$, in the former case, we have
\begin{align*}
\int^t_s\left| (t-s)^{-1} \Pi^*_{r, t}\bar\xi+\bar\eta\right|^2 \dif r
&\geq \int^{s+\eps (t-s)}_s\left| (t-s)^{-1} \Pi^*_{r, t}\bar\xi+\bar\eta\right|^2\dif r \\
&\geq \int^{s+ \eps (t-s)}_s\left(\tfrac{1}{2}\left| (t-s)^{-1} \Pi^*_{r, t}\bar\xi\right|^2-|\bar\eta|^2\right)\dif r\\
&\geq  \int^{s+ \eps (t-s)}_s\left( \frac{| \bar \xi|}{2 (t-s) \| \Pi_{r, t}^{-1} \| } \right)^2\dif r
-|\bar\eta|^2 \eps (t-s) \\
&\geq \left( \frac{1-\delta}{4 \kappa^2_0}-\delta\right) \eps (t-s),
\end{align*}
and in the later case,
\begin{align*}
\int^t_s\left| (t-s)^{-1} \Pi^*_{r, t}\bar\xi+\bar\eta\right|^2\dif r
&\geq   \int^t_{t-\eps(t-s)}\left(\tfrac{1}{2}|\bar\eta|^2-\left| (t-s)^{-1} \Pi^*_{r, t}\bar\xi \right|^2\right)\dif r\\
&\geq  \left( \tfrac{1}{2}\delta -\|U\|_\infty^2\eps^{2} \right) \eps (t-s).
\end{align*}
Combining the above two cases, by first choosing $\delta$ small enough and then $\eps$ small enough, one finds that for some
$c_2=c_2(\alpha,\kappa_0)>0$,
\begin{align}\label{HJ1}
\inf_{|\bar\xi|^2+|\bar\eta|^2=1}
\int^t_s\left|(t-s)^{-1} \Pi^*_{r, t}\bar\xi+\bar\eta\right|^2 \geq c_2 (t-s),
\end{align}
which together with \eqref{NB8} gives
\begin{align}
\int^t_s \psi^{\nu_1} \left(\sigma^*_{r, t}\Pi^*_{r,t}\xi+\sigma^*_{r, t} \eta\right)\dif s
\geq c_3 (t-s) \left( |((t-s) \xi,\eta)|^2 \wedge |((t-s) \xi,\eta)|^\alpha \right).
\label{NB2}
\end{align}
Hence by \eqref{NB9},
$$
\mE\e^{\mathrm{i}\<(\xi,\eta), K^{\nu_1}_{s,t}\>} \leq\exp \left( - c_3 (t-s) \left( |((t-s) \xi,\eta)|^2 \wedge |((t-s) \xi,\eta)|^\alpha \right)
\right).
$$
On the other hand, by \eqref{FG1}, one sees that $\psi^{\nu_1}$ is smooth and for any $k\in\mN$,
\begin{align}
|\nabla^k\psi^{\nu_1}(\xi)|\leq C(|\xi|^m+1),\ \xi\in\mR^d\label{NB1}
\end{align}
for some $m\in\mN$ and $C>0$.
Thus \eqref{NB9}, \eqref{NB2} and \eqref{NB1} in particular implies that
$$
(\xi,\eta)\mapsto \mE\e^{\mathrm{i}\<(\xi,\eta), K^{\nu_1}_{0, 1}\>}\in\cS(\mR^{2d}).
$$
Therefore, $K^{\nu_1}_{0, 1}$ has a smooth density $p_1^{\nu_1}(x,\v)\in\cS(\mR^{2d})$, which is given by the inverse Fourier transform
$$
p_1^{\nu_1}(x,\v)=\int_{\mR^{2d}}\e^{-\mathrm{i}(\<x,\xi\>+\<\v,\eta\>)}\mE\e^{\mathrm{i}\<(\xi,\eta), K^{\nu_1}_{0, 1} \>}\dif\xi\dif\eta.
$$
In particular, \eqref{LK33} holds.
\medskip\\
{\bf (iv)} Finally, we assume \eqref{Con2}, and make the following decomposition
$$
\nu_s=\nu^{(\alpha)}_1+\mu_s,
$$
where $\mu_s:=\nu_s-\nu^{(\alpha)}_1\in\mL^{sym}$.
Let $L^{\nu^{(\alpha)}}$ be an $\alpha$-stable process with the L\'evy measure $\nu^{(\alpha)}_1$, and let
$N_0(\dif t,\dif\v)$ be an independent Poisson random measure with intensity measure $\mu_t(\dif\v)\dif t$.
Let $L^\mu_{s,t}$ be defined as in \eqref{ES1} and $K^{\mu}_{s,t}=(X^{\mu}_{s,t},V^{\mu}_{s,t})$ be defined as in \eqref{KT0}  with
$L^{\mu}_{s,t}$ in place of $L^\nu_{s,t}$. Clearly,
\begin{align}
(L^\nu_{s,t})_{s\leq t}\stackrel{(d)}{=}(L^{\nu^{(\alpha)}_1}_{t-s}+L^{\mu}_{s,t})_{s\leq t}, \ \
K^\nu_{s,t}\stackrel{(d)}{=}K^{\nu^{(\alpha)}_1}_{s,t}+K^{\mu}_{s,t}.\label{LK1}
\end{align}
Thus, the distributional density of $K^\nu_{s,t}$ is given by
\begin{align}\label{EX1}
p^\nu_{s,t}(x,\v)=\mE p^{\nu^{(\alpha)}_1}_{s,t}\big(x-X^{\mu}_{s,t},\v-V^{\mu}_{s,t}\big).
\end{align}
As above, by \eqref{Mom} and $|\Pi_{s,t}|\leq \|U\|_\infty(t-s)$, we have
$$
\mE|X^{\mu}_{s,t}|^{q_1}\leq C(t-s)^{q_1+\frac{q_1}{\alpha}},\ \mE|V^{\mu}_{s,t}|^{q_2}\leq C(t-s)^{\frac{q_2}{\alpha}},
$$
and
$$
\mE\Big(|X^{\mu}_{s,t}|^{q_1}|V^{\mu}_{s,t}|^{q_2}\Big)\leq\left(\mE|X^{\mu}_{s,t}|^{q_1+q_2}\right)^{\frac{q_1}{q_1+q_2}}
\left(\mE|V^{\mu}_{s,t}|^{q_1+q_2}\right)^{\frac{q_2}{q_1+q_2}}\leq C(t-s)^{q_1+\frac{q_1+q_2}{\alpha}},
$$
which, together with \eqref{EX1} and what we have proved, gives \eqref{EV11}.
\end{proof}
\br
Let $p_{s,t}(x',\v'; x,\v)$ be the smooth density of $K_{s,t}(x,\v)=K_{s,t}+(x+\Pi_{s,t}\v, \v)$, which is given by
\begin{align}\label{Den}
p_{s,t}(x',\v'; x,\v)=p^{\nu}_{s,t}\big(x'-x-\Pi_{s,t}\v,\v'-\v\big).
\end{align}
For any $n_1,m_1,n_2,m_2\in\mN_0$, there is a constant $C=C(d,n_i,m_i,\kappa_0,\nu^{(\alpha)}_1,\alpha)>0$ such that for all $s<t$ and $x,\v\in\mR^d$,
\begin{align}
\int_{\mR^{2d}}|\nabla^{n_1}_{x'}\nabla^{m_1}_{\v'}\nabla^{n_2}_{x}\nabla^{m_2}_{\v} p_{s,t}(x',\v';x,\v)|
\dif x'\dif \v'\leq C(t-s)^{-((n_1+n_2)(1+\alpha)+m_1+m_2)/{\alpha}},\label{EV1}
\end{align}
which follows by the chain rule, $|\Pi_{s,t}|\leq\|U\|_\infty (t-s)$ and \eqref{EV11}.
\er
\bc
Under \eqref{Kapp} and \eqref{Con2},
for any $f\in \cB_b(\mR^{2d})$, $\cT_{s,t}f$ satisfies the following backward Kolmogorov's equation: for Lebesgue-almost all $s<t$ and all $x,\v\in\mR^d$,
\begin{align}
\p_s\cT_{s,t} f(x,\v)+\sK_s\cT_{s,t} f(x,\v)=0,\label{EQ0}
\end{align}
where $\sK_s$ is defined by \eqref{Eq218}.
\ec
\begin{proof}
First of all, as a consequence of (\ref{EV1}), we have for any $n,m\in\mN_0$,
\begin{align}
\|\nabla^n_x\nabla^m_\v\cT_{s,t} f\|_\infty\leq C(t-s)^{- (n(1+\alpha)+m)/{\alpha}}\|f\|_\infty,\ \ s<t.\label{EV6}
\end{align}
Thus, by Lebesgue's differentiable theorem,
it suffices to prove that for all $s\leq t_0<t$ and all $x,\v\in\mR^d$,
$$
\cT_{s,t} f(x,\v)=\cT_{t_0,t} f(x,\v)+\int^{t_0}_s\sK_r\cT_{r,t} f(x,\v)\dif r.
$$
Fix $t_1\in(t_0,t)$ and define $g(x,\v):=\cT_{t_1,t} f(x,\v)$. By \eqref{Sem}, we only need to show that  for all $s\leq t_0$ and all $x,\v\in\mR^d$,
$$
\cT_{s,t_1}g(x,\v)=\cT_{t_0,t_1} g(x,\v)+\int^{t_0}_s\sK_r\cT_{r,t_1}g(x,\v)\dif r.
$$
Since $g\in C^\infty_b(\mR^{2d})$ by \eqref{EV6}, it follows by \eqref{EQ00}.
\end{proof}

\bl\label{Le27}
Let $\beta,\gamma\in(0,2)$. Under \eqref{Kapp} and \eqref{Con2},
for any $\bar\sigma\in\mM^d$ and $\bar\nu\in\mL^{sym}$ with $\|\bar\sigma\|\leq\kappa_0$ and $\bar\nu\leq\nu^{(\gamma)}\in\mL^{sym, (\gamma)}$,
and for any $n,m\in\mN_0$,
there is a positive constant $C$ depending only on
$\kappa_0, \nu, n,m,d,\nu^{(\alpha)}_i,\nu^{(\gamma)}, \beta,\gamma$ such that for any $f\in C^2_b(\mR^{2d})$ and $t>s$,
\begin{align}
\big\|\nabla^n_x\nabla^m_\v\sL^{\bar\nu}_{\bar\sigma,\v}\cT_{s,t}\Delta^{\frac{\beta}{2}}_\v f\big\|_\infty
&\leq C(t-s)^{-(n(1+\alpha)+m+\beta+\gamma)/{\alpha}}\|f\|_\infty,\label{JH2}\\
\big\|\nabla^n_x\nabla^m_\v\sL^{\bar\nu}_{\bar\sigma,\v}\Delta^{\frac{\beta}{2}}_x\cT_{s,t} f\big\|_\infty
&\leq C(t-s)^{-(n(1+\alpha)+m+\beta+\gamma )/{\alpha}-\beta}\|f\|_\infty.\label{JH1}
\end{align}
Here we use the convention: $\sL^{0}_{\bar\sigma,\v}\equiv\mI$ the identity operator.
\el

\begin{proof}
Let $p_{s,t}(x',\v',x,\v)$ be given by \eqref{Den}. Notice that  by definition,
\begin{align*}
\nabla^n_x\nabla^m_\v\sL^{\bar\nu}_{\bar\sigma,\v}\cT_{s,t}\Delta^{\frac{\beta}{2}}_\v f(x,\v)
&=\int_{\mR^{2d}}\nabla^n_x\nabla^m_\v\sL^{\bar\nu}_{\bar\sigma,\v} p_{s,t}(x',\v',x,\v)\Delta^{\frac{\beta}{2}}_{\v'} f(x',\v')\dif x'\dif\v'\\
&=\int_{\mR^{2d}}\Delta^{\frac{\beta}{2}}_{\v'}\nabla^n_x\nabla^m_\v\sL^{\bar\nu}_{\bar\sigma,\v} p_{s,t}(x',\v',x,\v) f(x',\v')\dif x'\dif\v',
\end{align*}
and
\begin{align*}
\Delta^{\frac{\beta}{2}}_{\v'}\nabla^n_x\nabla^m_\v\sL^{\bar\nu}_{\bar\sigma,\v} p_{s,t}(x',\v',x,\v)
\!=\!\!\!\int_{\mR^d}\!\int_{\mR^d}\nabla^n_x\nabla^m_\v\delta^{(2)}_{\bar\v'}\delta^{(2)}_{\bar\sigma\bar\v}
p_{s,t}(x',\v',x,\v)\bar\nu(\dif\bar\v)\frac{\dif\bar \v'}{|\bar\v'|^{d+\beta}}.
\end{align*}
By using \eqref{EV311}, \eqref{EV31} and \eqref{EV1}, it is easy to see that for some $C>0$ independent of $x,\v,\bar\v,\bar\v'$,
\begin{align*}
&\int_{\mR^{2d}}\big|\nabla^n_x\nabla^m_\v\delta^{(2)}_{\bar\v'}\delta^{(2)}_{\bar\sigma\bar\v} p_{s,t}(x',\v',x,\v)\big|\dif x'\dif\v'
\leq C(t-s)^{-\frac{n(1+\alpha)+m}{\alpha}}\\
&\qquad\times \Big(\big((t-s)^{-\frac{4}{\alpha}}|\bar\v'|^2|\bar\v|^2\big)\wedge
\big((t-s)^{-\frac{2}{\alpha}}(|\bar\v'|^2\wedge |\bar\v|^2)\big)\wedge1\Big).
\end{align*}
Hence,
\begin{align*}
&\|\nabla^n_x\nabla^m_\v\cL^{\bar\nu}_{\bar\sigma,\v}\cT_{s,t}\Delta^{\frac{\beta}{2}}_\v f\|_\infty
\leq C(t-s)^{-\frac{m+n(1+\alpha)}{\alpha}}\|f\|_\infty\\
&\times\int_{\mR^d}\!\!\int_{\mR^d}\Big(\big((t-s)^{-\frac{4}{\alpha}}|\bar\v'|^2|\bar\v|^2\big)\wedge
\big((t-s)^{-\frac{2}{\alpha}}(|\bar\v'|^2\wedge |\bar\v|^2)\big)\wedge1\Big)\nu^{(\gamma)}(\dif\bar\v)\frac{\dif\bar \v'}{|\bar\v'|^{d+\beta}}.
\end{align*}
If we calculate the double integral in the following four regions separately,
\begin{align*}
&\Big\{\bar\v|\leq(t-s)^{\frac{1}{\alpha}},|\bar\v'|\leq(t-s)^{\frac{1}{\alpha}}\Big\}
\cup\Big\{\bar\v|\leq(t-s)^{\frac{1}{\alpha}},|\bar\v'|>(t-s)^{\frac{1}{\alpha}}\Big\}\\
&\cup\Big\{\bar\v|>(t-s)^{\frac{1}{\alpha}},|\bar\v'|\leq(t-s)^{\frac{1}{\alpha}}\Big\}
\cup\Big\{\bar\v|>(t-s)^{\frac{1}{\alpha}},|\bar\v'|>(t-s)^{\frac{1}{\alpha}}\Big\},
\end{align*}
then we obtain \eqref{JH2}. Similarly, one can show \eqref{JH1}.
\end{proof}

\subsection{Fefferman-Stein's theorem}
In this subsection we recall the classical Fefferman-Stein's theorem. First of all, we introduce a family of ``balls'' looking like a ``parallelepiped'' in $\mR^{1+2d}$,
as seen below, which is natural for treating the kinetic operator.
More precisely, fixing $\alpha\in(0,2)$, and for any $r>0$ and point $(t_0,x_0,\v_0)\in\mR^{1+2d}$, we define
\begin{align}
Q_r(t_0,x_0,\v_0):=\Big\{(t,x,\v):\ &t\in B_{r^\alpha}(t_0),x\in B_{r^{1+\alpha}}\big(x_0+\Pi_{t_0,t}\v_0\big), \v\in B_r(\v_0)\Big\},\label{Ball}
\end{align}
where $\Pi_{t_0,t}:=\int^t_{t_0}U_r\dif r$ and $B_r(\v_0)$ is the Euclidean ball with radius $r$ and center $\v_0$.
The set of all such balls is denoted by $\mQ^{(\alpha)}$.
For $f\in L^1_{loc}(\mR^{1+2d})$, we define the Hardy-Littlewood maximal function by
$$
\cM f(t,x,\v):=\sup_{r>0}\fint_{Q_r(t,x,\v)}|f(t',x',\v')|\dif \v'\dif x'\dif t',
$$
and the sharp function by
$$
\cM^\sharp f(t,x,\v):=\sup_{r>0}\fint_{Q_r(t,x,\v)}|f(t',x',\v')-f_{Q_r(t,x,\v)}|\dif \v'\dif x'\dif t',
$$
where for a $Q\in\mQ^{(\alpha)}$, $|Q|$ denotes the Lebesgue measure of $Q$ and
$$
f_Q:=\fint_Q f(t',x',\v')\dif \v'\dif x'\dif t'=\frac{1}{|Q|}\int_Q f(t',x',\v')\dif \v'\dif x'\dif t'.
$$
One says that a function $f\in BMO(\mR^{1+2d})$ if $\cM^\sharp f\in L^\infty(\mR^{1+2d})$. Clearly, $f\in BMO(\mR^{1+2d})$ if and only if
there exists a constant $C>0$ such that for any $Q\in\mQ^{(\alpha)}$, and for some $a_Q\in\mR$,
$$
\fint_Q |f(t',x',\v')-a_Q|\dif \v'\dif x'\dif t'\leq C.
$$

We have the following simple property about $Q_r\in\mQ^{(\alpha)}$.
\bp\label{Pr24}
Let $c_1:=3^\frac{1}{\alpha}\vee 3\vee (3+4\|U\|_\infty)^{\frac{1}{1+\alpha}}$ and $c_2:=c_1^{1+(2+\alpha)d}$. We have
\begin{enumerate}[(i)]
\item If $Q_r(t_0,x_0,\v_0)\cap Q_r(t_0',x_0',\v_0')\not=\emptyset$, then
\begin{align}
Q_r(t_0,x_0,\v_0)\subset Q_{c_1 r}(t_0',x_0',\v_0').\label{EL1}
\end{align}
\item $|Q_{c_1 r}(t_0,x_0,\v_0)|\leq c_2|Q_{r}(t_0,x_0,\v_0)|$.
\end{enumerate}
\ep
\begin{proof}
(i) By the assumption, we have
$$
|t_0-t_0'|\leq 2r^{\alpha}, \ |\v_0-\v_0'|\leq 2r,
$$
and for some $t'\in B_{r^\alpha}(t_0)\cap B_{r^\alpha}(t_0')$,
$$
\big|x_0-\Pi_{t_0,t'}\v_0-\big(x_0'-\Pi_{t'_0,t'}\v'_0\big)\big|\leq 2r^{1+\alpha}.
$$
Thus, for any $(t,x,\v)\in Q_r(t_0,x_0,\v_0)$, we have
$$
|t-t_0'|\leq 3r^\alpha,\ |\v-\v_0'|\leq 3r
$$
and
\begin{align*}
\big|x-\big(x_0'-\Pi_{t_0',t}\v_0'\big)\big|&\leq \big|x-\big(x_0-\Pi_{t_0,t}\v_0\big)\big|+\big|x_0-\Pi_{t_0,t}\v_0-\big(x_0'-\Pi_{t_0',t}\v_0'\big)\big|\\
&\leq r^{1+\alpha}+\big|x_0-\Pi_{t_0,t'}\v_0-\big(x_0'-\Pi_{t'_0,t'}\v_0'\big)\big|+\big|\Pi_{t,t'}(\v_0-\v_0')\big|\\
&\leq r^{1+\alpha}+2r^{1+\alpha}+4\|U\|_\infty r^{1+\alpha}=(3+4\|U\|_\infty)r^{1+\alpha}.
\end{align*}
From these, we immediately obtain (\ref{EL1}).
\\
\\
(ii) It follows by noticing that $|Q_{r}(t_0,x_0,\v_0)|=c_3 r^{1+(2+\alpha)d}$  for some $c_3=c_3(d)$.
\end{proof}
\br
By Proposition \ref{Pr24} and \cite[Theorem 1, p.13]{St}, for any $p\in(1,\infty]$, there is a constant $C>0$ such that for any $f\in L^p(\mR^{1+2d})$,
\begin{align}
\|\cM f\|_p\leq C\|f\|_p.\label{EN1}
\end{align}
\er
We need the following version of Fefferman-Stein's theorem, whose proof is given in Appendix 4.2.
\bt\label{Th26}
(Fefferman-Stein's theorem) For any $p\in(1,\infty)$, there exists
a constant $C=C(p,d,\alpha)>0$ such that for all $f\in L^p(\mR^{1+2d})$,
\begin{align}
\|f\|_p\leq C\|\cM^\sharp f\|_p.\label{EU4}
\end{align}
\et
Using this theorem, we have
\bt\label{Th2}
For $q\in(1,\infty)$, let $\sP$ be a bounded linear operator from $L^q(\mR^{1+2d})$ to $L^q(\mR^{1+2d})$
and also from $L^\infty(\mR^{1+2d})$ to $BMO(\mR^{1+2d})$. Then for any $p\in[q,\infty)$ and $f\in L^p(\mR^{1+2d})$,
$$
\|\sP f\|_p\leq C\|f\|_p,
$$
where the constant $C$ depends only on $p,q$ and the norms of $\|\sP\|_{L^q\to L^q}$ and $\|\sP\|_{L^\infty\to BMO}$.
\et
\begin{proof}
Noticing that by the assumptions,
$$
\|\cM^\sharp(\sP f)\|_q\leq 2\|\cM(\sP f)\|_q\stackrel{(\ref{EN1})}{\leq} C\|\sP f\|_q\leq C\|\sP\|_{L^q\to L^q}\|f\|_q
$$
and
$$
\|\cM^\sharp(\sP f)\|_\infty\leq \|\sP\|_{L^\infty\to BMO}\|f\|_\infty,
$$
by the classical Marcinkiewicz's interpolation theorem (cf. \cite{St}), we have for any $p\in[q,\infty)$,
$$
\|\cM^\sharp(\sP f)\|_p\leq C\|f\|_p,
$$
which together with (\ref{EU4}) gives the desired estimate.
\end{proof}

\section{$L^p$-maximal regularity of nonlocal kinetic Fokker-Planck equations}

For $\lambda>0$, we consider the following linear equation:
\begin{align}
\p_s u+(\sK_s-\lambda) u+f=0,\label{EQ101}
\end{align}
where $\sK_s$ is defined by \eqref{Eq218}.
We first introduce the following notion.
\bd\label{Def31}
For given $f\in L^1_{loc}(\mR^{1+2d})$, a function $u\in C(\mR; L^1_{loc}(\mR^{2d}))$ is called a weak solution of equation \eqref{EQ101}
if for all $s\leq T$ and any $\phi\in C^\infty_c(\mR^{2d})$,
\begin{equation}\label{Weak}
\<u(s),\phi\>=\<u(T),\phi\>+\int^T_s\<u(t),(\sK^*_t-\lambda) \phi\>\dif t+\int^T_s\<f(t),\phi\>\dif t,
\end{equation}
where $\<u,\phi\>:=\int_{\mR^{2d}}u(x,\v)\phi(x,\v)\dif x\dif \v$ and $\sK^*_t:=\sL^{\nu_t}_{\sigma_t,\v}-U_t\v\cdot\nabla_x$ is the adjoint operator of $\sK_t$.
\ed
We need the following simple result.
\bp\label{Pr32}
Given $p\in[1,\infty]$ and $f\in L^p(\mR^{1+2d})$, the unique weak solution of equation \eqref{EQ101} with $u\in C(\mR; L^p(\mR^{2d}))$
 and $\lim_{T_n\to \infty} u(T_n)=0$ weakly for some deterministic sequence $T_n\to \infty$
 is given by
\begin{align}
u(s,x,\v)=\int^{\infty}_s\e^{\lambda(s-t)}\cT_{s,t} f(t,x,\v)\dif t,\label{TST1}
\end{align}
where $\cT_{s,t}f$ is defined by \eqref{TST}.
\ep
\begin{proof}
Let $\varrho:\mR^{2d}\to[0,\infty)$ be a smooth function with compact support and $\int\varrho=1$. For $\eps>0$, define
$$
\varrho_\eps(x,\v):=\eps^{-3d}\varrho(\eps^{-1}x, \eps^{-2}\v),\ \ f_\eps(t,x,\v):=f(t)*\varrho(x,\v),
$$
where $*$ denotes the convolution, and
$$
u_\eps(s,x,\v):=\int^{\infty}_s\e^{\lambda(s-t)}\cT_{s,t} f_\eps(t,x,\v)\dif t.
$$
Since $f_\eps\in L^p(\mR; C^\infty_b(\mR^{2d}))$,  we have by \eqref{EQ00},
$$
\p_s u_\eps+(\sK_s-\lambda) u_\eps+f_\eps=0.
$$
In particular, for all $s\leq T$ and $\phi\in C^\infty_c(\mR^{2d})$,
\begin{align*}
\<u_\eps(s),\phi\>=\<u_\eps(T),\phi\>+\int^T_s\<u_\eps(t),(\sK^*_t-\lambda) \phi\>\dif t+\int^T_s\<f_\eps(t),\phi\>\dif t.
\end{align*}
By taking limits $\eps\to 0$ and the dominated convergence theorem, one sees that $u$ is a weak solution of equation \eqref{EQ101}.
Moreover, we also have $u\in C(\mR; L^p(\mR^{2d}))$ and $\lim_{T\to\infty} u(T)=0$
weakly.

On the other hand, let $u$ be a weak solution of \eqref{EQ101}. In \eqref{Weak}, taking $\phi=\varrho_\eps(x-\cdot,\v-\cdot)$
and setting $u_\eps:=u*\varrho_\eps$, $f_\eps:=f*\varrho_\eps$, one has
\begin{align}
u_\eps(s)=u_\eps(T)+\int^T_s(\sK_t-\lambda)u_\eps(t)\dif t+\int^T_s(f_\eps+g_\eps)(t)\dif t,\label{BN9}
\end{align}
where
\begin{align*}
g_\eps(t,x,\v)&:=\int_{\mR^{2d}} u(t,x',\v') U_t(\v'-\v)\cdot\nabla_x\varrho_\eps(x-x',\v-\v')\dif x'\dif \v'\\
&=\int_{\mR^{2d}} u(t,x-x',\v-\v')U_t\v'\cdot\nabla_x\varrho_\eps(x',\v')\dif x'\dif \v'.
\end{align*}
Since $u_\eps\in C(\mR; C^\infty_b(\mR^{2d}))$ and
$\lim_{n\to \infty} u_\eps (T_n)=0$,
the unique solution of \eqref{BN9} is given by
\begin{align}
u_\eps(s,x,\v)=\int^{\infty}_s\e^{\lambda(s-t)}\cT_{s,t} (f_\eps+g_\eps)(t,x,\v)\dif t.\label{TST2}
\end{align}
Notice that by the definition of $\varrho_\eps$,
$$
\|g_\eps(t)\|_p\leq \|u(t)\|_p\int_{\mR^{2d}}|U_t\v'|\cdot|\nabla_x\varrho_\eps(x',\v')|\dif x'\dif \v'\leq C\eps\|u(t)\|_p\to 0.
$$
By taking limits $\eps\to 0$ for both sides of \eqref{TST2}, we obtain \eqref{TST1}.
\end{proof}

Now we can present our main result of this paper.

\bt\label{Main}
Under  \eqref{Kapp} and \eqref{Con2}, for any $p\in(1,\infty)$,
there exists a positive constant $C=C(\kappa_0, p,d,\nu^{(\alpha)}_i,\alpha)$ such that for all $\lambda> 0$ and $f\in L^p(\mR^{1+2d})$,
\begin{align}
\big\|\Delta_x^{\frac{\alpha}{2(1+\alpha)}}u^\lambda\big\|_p+\big\|\Delta_\v^{\frac{\alpha}{2}}u^\lambda\big\|_p\leq C\|f\|_p,\label{BN77}
\end{align}
where $u^\lambda(s,x,\v):=\int^\infty_s\e^{\lambda(s-t)}\cT_{s,t} f_t(x, \v)\dif t$ is the unique weak solution of equation \eqref{EQ101}.
\et

When $p=2$ and $U_s$ is independent of $s$, estimate \eqref{BN77} was proved in \cite{Al}.
The proof of Theorem \ref{Main} for $p=2$ will be given in Appendix 4.3, which is new and more elementary
even for the time-independent case considered in \cite{Al}.

\subsection{Proof of Theorem \ref{Main} for $p\in(2,\infty)$}
We  introduce the following two operators:
\begin{align*}
\sP_1f:=\sP^{\nu,\sigma, U}_1 f(s,x,\v)&:=\Delta_x^{\frac{\alpha}{2(1+\alpha)}}\int^\infty_s\e^{\lambda(s-t)}\cT^{\nu,\sigma, U}_{s,t} f(t,x,\v)\dif t,\\
\ \ \sP_2f:=\sP^{\nu,\sigma, U}_2 f(s,x,\v)&:=\Delta_\v^{\frac{\alpha}{2}}\int^\infty_s\e^{\lambda(s-t)}\cT^{\nu,\sigma, U}_{s,t} f(t,x,\v)\dif t.
\end{align*}
By Theorem \ref{Th2} and \eqref{BN77} for $p=2$,
our main task is to show that $\sP_1$ and $\sP_2$ are bounded linear operators from $L^\infty(\mR^{1+2d})$ to $BMO$.
More precisely,  we want to prove that for any $f\in L^\infty(\mR^{1+2d})$ with $\|f\|_\infty\leq 1$,
and any $Q=Q_r(t_0,x_0,\v_0)\in\mQ^{(\alpha)}$,
\begin{align}
\fint_{Q}|\sP^{\nu,\sigma,U}_i f(s,x,\v)-a^{Q}_i|^2 \leq C,\label{ET0}
\end{align}
where $a^{Q}_i$ is a constant depending on $Q$ and $f$, and $C$ only depends on $\kappa_0, p,d,\nu^{(\alpha)}_i,\alpha$.
\bl\label{Le34}
(Scaling Property) For any $Q=Q_r(t_0,x_0,\v_0)\in\mQ^{(\alpha)}$ and $i=1,2$, we have
\begin{align}
\fint_{Q_r(t_0,x_0,\v_0)}\big|\sP^{\nu,\sigma,U}_if(s,x,\v)-a\big|^2=
\fint_{Q_1(0)}\big|\sP^{\tilde\nu,\tilde\sigma,\tilde U}_i\tilde f(s,x,\v)-a\big|^2,\label{SCL}
\end{align}
where $a\in\mR$, $\tilde\nu_s:=\nu_{r^\alpha s+t_0}$, $\tilde \sigma_s:=\sigma_{r^\alpha s+t_0}$, $\tilde U_s:=U_{r^\alpha s+t_0}$ and
$$
\tilde f(t,x,\v):=f\big(r^\alpha t+t_0, r^{1+\alpha}x+x_0+\Pi_{t_0,r^\alpha t+t_0} \v_0, r\v+\v_0\big).
$$
\el
\begin{proof}
Let us write
$$
u(s,x,\v):=\int^\infty_s\e^{\lambda(s-t)}\cT^{\nu,\sigma, U}_{s,t} f(t,x,\v)\dif t
$$
and
$$
\tilde u(s,x,\v):=r^{-\alpha} u\big(r^\alpha s+t_0,r^{1+\alpha} x+x_0+\Pi_{t_0,r^\alpha s+t_0} \v_0, r\v+\v_0\big),
$$
where $\Pi_{t_0,r^\alpha t+t_0}=\int^{r^\alpha t+t_0}_{t_0}U_{r'}\dif r'$. By the change of variables, we have
\begin{align*}
\fint_{Q_r(t_0,x_0,\v_0)}\big|\Delta_x^{\frac{\alpha}{2(1+\alpha)}}u(s,x,\v)-a\big|^2
=\fint_{Q_1(0)}\big|\Delta_x^{\frac{\alpha}{2(1+\alpha)}}\tilde u(s,x,\v)-a\big|^2.
\end{align*}
On the other hand, by Proposition \ref{Pr32}, one sees that
$$
\tilde u(s,x,\v)=\int^\infty_s\e^{\lambda(s-t)}\cT^{\tilde \nu,\tilde \sigma, \tilde U}_{s,t}\tilde f(t,x,\v)\dif t.
$$
Thus, we obtain \eqref{SCL} for $i=1$. Similarly, \eqref{SCL} holds for $i=2$.
\end{proof}

Below we split $\sP_if=\sP_{i1}f+\sP_{i2}f$, $i=1,2$, where
\begin{align*}
\sP_{11} f&:=\Delta_x^{\frac{\alpha}{2(1+\alpha)}}\int^2_\cdot\e^{\lambda(\cdot-t)}\cT_{\cdot,t}f(t)\dif t,\ \ \
\sP_{21} f:=\Delta_\v^{\frac{\alpha}{2}}\int^2_\cdot\e^{\lambda(\cdot-t)}\cT_{\cdot,t}f(t)\dif t,\\
\sP_{12} f&:=\Delta_x^{\frac{\alpha}{2(1+\alpha)}}\int^\infty_2\e^{\lambda(\cdot-t)}\cT_{\cdot,t} f(t)\dif t,\ \ \
\sP_{22} f:=\Delta_\v^{\frac{\alpha}{2}}\int^\infty_2\e^{\lambda(\cdot-t)}\cT_{\cdot,t} f(t)\dif t.
\end{align*}
First of all, we treat $\sP_{11}f, \sP_{21}f$.
\bl\label{Le35}
Under  \eqref{Kapp} and \eqref{Con2}, there is a constant $C>0$ depending only on $\kappa_0, p,d,\nu^{(\alpha)}_i,\alpha$
such that for all $f\in L^\infty(\mR^{1+2d})$ with $\|f\|_\infty\leq 1$,
\begin{align}
\int_{Q_1(0)}|\sP_{i1} f(s,x,\v)|^2\leq C,\ \ i=1,2.\label{BN4}
\end{align}
\el
\begin{proof}
For $s\in[-1,1]$, let
$$
u(s,x,\v):=\int^2_s\e^{\lambda(s-t)}\cT_{s,t}f (t,x,\v)\dif t=\int^\infty_s\e^{\lambda(s-t)}\cT_{s,t}((1_{[-1,2]}f) (t))(x,\v)\dif t.
$$
Since $\|f\|_\infty\leq 1$, we have
\begin{align}
\|u(s)\|_\infty\leq 3,\ \ s\in[-1,1].\label{BOU}
\end{align}
By \eqref{EV311}, \eqref{EV31} and \eqref{EV6},  we have for any $t>s$,
\begin{align*}
\|\delta^{(1)}_x\cT_{s,t}f (t)\|_\infty
&\preceq (\|\nabla_x\cT_{s,t}f (t)\|_\infty|x|)\wedge\|\cT_{s,t}f (t)\|_\infty\preceq ((t-s)^{-\frac{1+\alpha}{\alpha}}|x|)\wedge 1,\\
\|\delta^{(1)}_\v\cT_{s,t}f(t)\|_\infty
&\preceq (\|\nabla_\v\cT_{s,t}f (t)\|_\infty|\v|)\wedge\|\cT_{s,t}f (t)\|_\infty\preceq ((t-s)^{-\frac{1}{\alpha}}|\v|)\wedge 1.
\end{align*}
Since $a\wedge 1\leq a^\gamma$ for any $a>0$ and $\gamma\in[0,1]$,  we have for any $\gamma_1\in(0,\alpha/(1+\alpha))$,
\begin{align}
\|\delta^{(1)}_xu(s)\|_\infty\preceq |x|^{\gamma_1}\int^2_s\e^{\lambda(s-t)}(t-s)^{-(1+\alpha)\gamma_1/\alpha}\dif t\leq C |x|^{\gamma_1},\label{BN2}
\end{align}
and for any $\gamma_2\in(0,\alpha\wedge 1)$,
\begin{align}
\|\delta^{(1)}_\v u(s)\|_\infty&\preceq |\v|^{\gamma_2}\int^2_s\e^{\lambda(s-t)}(t-s)^{-\gamma_2/\alpha}\dif t\leq C|\v|^{\gamma_2},\label{BN11}
\end{align}
where $C>0$ is independent of $\lambda>0$.

Let $\varphi$ be a nonnegative smooth cutoff function in $\mR^{2d}$ with $\varphi(x,\v)=1$ for $|(x,\v)|\leq 4$ and  $\varphi(x)=0$ for $|(x,\v)|>8$.
By Definition \ref{Def31}, it is easy to see that
$u\varphi$ is a weak solution of equation \eqref{EQ101} with $f$ replacing by
$$
g_\varphi=\left(f\varphi+\sK_s\varphi\, u+\int_{\mR^d}\delta^{(1)}_\v u\,\delta^{(1)}_\v\varphi\,\nu_s(\dif\v)\right)1_{[-1,2]}(s).
$$
Noticing that by \eqref{BOU} and \eqref{BN11},
$$
\|g_\varphi\|_2\leq C_\varphi,
$$
and by Proposition \ref{Pr32}, we have
$$
(u\varphi)(s,x,\v)=\int^\infty_s\e^{\lambda(s-t)}\cT_{s,t}g_\varphi(t,x,\v)\dif t,
$$
which implies by \eqref{BN77} for $p=2$ that
\begin{align}
\big\|\Delta_x^{\frac{\alpha}{2(1+\alpha)}}(u\varphi)\big\|_2+\big\|\Delta_\v^{\frac{\alpha}{2}}(u\varphi)\big\|_2\leq C\|g_\varphi\|_2\leq C.\label{BN5}
\end{align}
By the definition of $\sP_{11}$ and \eqref{EV441}, \eqref{BN2},  \eqref{BN5}, we have
\begin{align*}
&\int_{Q_1(0)}|\sP_{11} f|^2=\int_{Q_1(0)}|\Delta_x^{\frac{\alpha}{2(1+\alpha)}}u|^2\leq\int_{\mR^{1+2d}}|(\Delta_x^{\frac{\alpha}{2(1+\alpha)}}u)\varphi|^2\\
&\quad\leq 2\big\|\Delta_x^{\frac{\alpha}{2(1+\alpha)}}(u\varphi)\big\|_2^2+C_\varphi\Big(\sup_{x}\|\delta^{(1)}_x u\|_\infty/|x|^{\gamma_1}+\|u\|_\infty\Big)\leq C,
\end{align*}
and by \eqref{EV441}, \eqref{BN11} and \eqref{BN5},
\begin{align*}
&\int_{Q_1(0)}|\sP_{21} f|^2=\int_{Q_1(0)}|\Delta_\v^{\frac{\alpha}{2}}u|^2\leq\int_{\mR^{1+2d}}|(\Delta_\v^{\frac{\alpha}{2}}u)\varphi|^2\\
&\quad\leq 2\big\|\Delta_\v^{\frac{\alpha}{2}}(u\varphi)\big\|_2^2+C_\varphi\Big(\sup_{x}\|\delta^{(1)}_\v u\|_\infty/|\v|^{\gamma_2}+\|u\|_\infty\Big)\leq C.
\end{align*}
The proof is complete.
\end{proof}

To treat $\sP_{12}f, \sP_{22}f$, we need the following estimate.
\bl
Under  \eqref{Kapp} and \eqref{Con2},
there is a constant $C>0$ depending only on $\kappa_0, p,d,\nu^{(\alpha)}_i,\alpha$
such that for all $f\in L^\infty(\mR^{1+2d})$ with $\|f\|_\infty\leq 1$ and all $s\in[-1,1]$,
\begin{align}
&\int^\infty_2\Big|\Delta_x^{\frac{\alpha}{2(1+\alpha)}}\cT_{s,t}f(t,0,0)
-\Delta_x^{\frac{\alpha}{2(1+\alpha)}}\cT_{0,t}f(t,0,0)\Big|\dif t\leq C,\label{NB6'}\\
&\quad\int^\infty_2\Big|\Delta_\v^{\frac{\alpha}{2}}\cT_{s,t}f(t,0,0)
-\Delta_\v^{\frac{\alpha}{2}}\cT_{0,t}f(t,0,0)\Big|\dif t\leq C.\label{NB6}
\end{align}
\el
\begin{proof}
First of all, by \eqref{JH1} with $\gamma=\alpha$ and $\beta=\frac{\alpha}{1+\alpha}$, we have for all $s\in[-1,1]$,
\begin{align*}
&\int^\infty_2\big|\Delta_x^{\frac{\alpha}{2(1+\alpha)}}\cT_{s,t}f(t,0,0)
-\Delta_x^{\frac{\alpha}{2(1+\alpha)}}\cT_{0,t}f(t,0,0)\big|\dif t\\
&\quad\leq\int^\infty_2\!\!\!\int^s_0\big|\p_r\Delta_x^{\frac{\alpha}{2(1+\alpha)}}\cT_{r,t}f(t,0,0)\big|\dif r\dif t\\
&\quad=\int^\infty_2\!\!\!\int^s_0\big|\sL^{\nu_r}_{\sigma_r, \v}\Delta_x^{\frac{\alpha}{2(1+\alpha)}}\cT_{r,t}f(t,0,0)\big|\dif r\dif t\\
&\quad\preceq \int^\infty_2\!\!\int^s_0(t-r)^{-2}\dif r\dif t\preceq 1,
\end{align*}
which give \eqref{NB6'}.

Next we deal with \eqref{NB6}. Let $\chi$ be a smooth cutoff function with $\chi(s)=1$ for $s\in[0,1]$ and $\chi(s)=0$ for $s>3$.
Fix $\gamma\in(1,1+\tfrac{\alpha}{2-\alpha})$ and define
$$
h_t(\v):=\chi(|\v|/t^{\gamma/\alpha}),\ t>0, \ \v\in\mR^d.
$$
By definition, we have
\begin{align*}
\Delta_\v^{\frac{\alpha}{2}}\cT_{s,t}f(t,0,0)
=\int_{\mR^{d}}\delta^{(2)}_\v\cT_{s,t}f(t,0,0)\frac{\dif\v}{|\v|^{d+\alpha}}=I_1(s,t)+I_2(s,t),
\end{align*}
where
\begin{align*}
I_1(s,t)&:=\int_{\mR^d}\delta^{(2)}_\v \cT_{s,t}f(t,0,0)(1-h_{t-s}(\v))\frac{\dif\v}{|\v|^{d+\alpha}},\\
I_2(s,t)&:=\int_{\mR^d}\delta^{(2)}_\v \cT_{s,t}f(t,0,0)h_{t-s}(\v)\frac{\dif\v}{|\v|^{d+\alpha}}.
\end{align*}
Thus, we can write
\begin{align}
&\int^\infty_2\Big|\Delta_\v^{\frac{\alpha}{2}}\cT_{s,t}f(t,0,0)
-\Delta_\v^{\frac{\alpha}{2}}\cT_{0,t}f(t,0,0)\Big|\dif t\no\\
&\quad\leq \int^\infty_2|I_1(s,t)-I_1(0,t)|\dif t+\int^\infty_2|I_2(s,t)-I_2(0,t)|\dif t.\label{EV7}
\end{align}
In view of $\gamma>1$, we have for all $s\in[-1,1]$,
\begin{align}
\begin{split}
\int^\infty_2|I_1(s,t)|\dif t&\preceq \int^\infty_2\left(\int_{\mR^d}|(1-h_{t-s}(\v))|\frac{\dif\v}{|\v|^{d+\alpha}}\right)\dif t\\
&\preceq \int^\infty_2\left(\int_{|\v|>(t-s)^{\gamma/\alpha}}\frac{\dif\v}{|\v|^{d+\alpha}}\right)\dif t
\preceq\int^\infty_2(t-s)^{-\gamma}\dif t\preceq 1.
\end{split}\label{EV77}
\end{align}
On the other hand, let us write
\begin{align*}
\int^\infty_2|I_2(s,t)-I_2(0,t)|\dif t\leq\int^\infty_2\!\!\!\int^s_0|\p_r I_2(r,t)|\dif r\dif t\leq J_{1}+J_{2},
\end{align*}
where
\begin{align*}
J_{1}:=\int^\infty_2\!\!\!\int^s_0\!\!\!\int_{\mR^d}|\delta^{(2)}_\v\p_r\cT_{r,t}f(t,0,0)h_{t-r}(\v)|\frac{\dif\v}{|\v|^{d+\alpha}}\dif r\dif t,\\
J_{2}:=\int^\infty_2\!\!\!\int^s_0\!\!\!\int_{\mR^d}|\delta^{(2)}_\v\cT_{r,t}f(t,0,0)\p_rh_{t-r}(\v)|\frac{\dif\v}{|\v|^{d+\alpha}}\dif r\dif t.
\end{align*}
Recalling definition \eqref{NB5}, by \eqref{EQ0} and (\ref{JH2}),  we have
\begin{align*}
&\quad|\delta^{(2)}_\v\p_r\cT_{r,t} f(t,0,0)|\\
&\leq|\delta^{(2)}_\v\sL^{\nu_r}_{\sigma_r, \v}\cT_{r,t} f(t,0,0)|+|(U_r\v\cdot \nabla_x)(\cT_{r,t} f(t,0,\v)-\cT_{r,t} f(t,0,-\v))|\\
&\leq 2\|\nabla^2_\v\sL^{\nu_r}_{\sigma_r,\v}\cT_{r,t}f(t)\|_\infty|\v|^2+2\|U\|_\infty|\v|^2\|\nabla_\v\nabla_x\cT_{r,t} f(t)\|_\infty
\preceq|\v|^2(t-r)^{-1-\frac{2}{\alpha}}.
\end{align*}
By Fubini's theorem, we have for all $s\in[-1,1]$,
\begin{align}
J_{1}&\preceq\int^\infty_2\!\!\!\int^{s}_0\left(\int_{\mR^d}
(t-r)^{-1-\frac{2}{\alpha}}h_{t-r}(\v)\frac{\dif\v}{|\v|^{d+\alpha-2}}\right)\dif r\dif t\no\\
&=\int^s_0\!\!\!\int^\infty_{2-r}\left(\int_{\mR^d}t^{-1-\frac{2}{\alpha}}h_{t}(\v)\frac{\dif\v}{|\v|^{d+\alpha-2}}\right)\dif t\dif r\no\\
&\preceq\int^\infty_1t^{-1-\frac{2}{\alpha}}t^{\frac{2\gamma}{\alpha}-\gamma}\dif t
\preceq 1\mbox{ since $\gamma\in(1,1+\tfrac{\alpha}{2-\alpha})$}.\label{EV8}
\end{align}
For $J_{2}$,  noticing that
$$
|\p_th_t(\v)|=\tfrac{\gamma}{\alpha}|\v|t^{-\gamma/\alpha-1}|\chi'(|\v|/t^{\gamma/\alpha})|
\leq \tfrac{\gamma}{\alpha} t^{-1}\|\chi'\|_\infty1_{\{t^{\gamma/\alpha}<|\v|<3t^{\gamma/\alpha}\}},
$$
we also have
\begin{align}
J_{2}\leq\int^\infty_1\left(\int_{\mR^d}|\p_th_t(\v)|\frac{\dif\v}{|\v|^{d+\alpha}}\right)\dif t
\preceq\int^\infty_1t^{-1-\gamma}\dif t\preceq1.\label{EV9}
\end{align}
Combining (\ref{EV7}), \eqref{EV77}, (\ref{EV8}) and (\ref{EV9}), we obtain \eqref{NB6}.
\end{proof}

Now, we treat $\sP_{12}f, \sP_{22}f$ as follows.
\bl\label{Le37}
Under  \eqref{Kapp} and \eqref{Con2}, there is a constant $C>0$ depending only on $\kappa_0, p,d,\nu^{(\alpha)}_i,\alpha$
such that for all $f\in L^\infty(\mR^{1+2d})$ with $\|f\|_\infty\leq 1$,
\begin{align}\label{EP1}
\int_{Q_1(0)}|\sP_{i2}f(s,x,\v)-\sP_{i2}f(0,0,0)|^2\leq C,\ \ i=1,2.
\end{align}
\el
\begin{proof}
For $i=2$, by definition, we have
\begin{align*}
&|\sP_{22}f(s,x,\v)-\sP_{22}f(0,0,0)|\leq\int^\infty_2|\e^{\lambda(s-t)}-\e^{-\lambda t}|\, \|\Delta_\v^{\frac{\alpha}{2}}\cT_{s,t}f\|_\infty\dif t\\
&\qquad\qquad+\int^\infty_2\e^{-\lambda t}|\Delta_\v^{\frac{\alpha}{2}}\cT_{s,t}f(x,\v)-\Delta_\v^{\frac{\alpha}{2}}\cT_{s,t}f(0,0)|\dif t\\
&\qquad\qquad+\int^\infty_2\e^{-\lambda t}|\Delta_\v^{\frac{\alpha}{2}}\cT_{s,t}f(0,0)-\Delta_\v^{\frac{\alpha}{2}}\cT_{0,t}f(0,0)|\dif t\\
&\qquad\qquad=:I_1(s)+I_2(s,x,\v)+I_3(s).
\end{align*}
Noticing that by Lemma \ref{Le27},
\begin{align*}
&\|\Delta_\v^{\frac{\alpha}{2}} \cT_{s,t}f(t)\|_\infty\leq C(t-s)^{-1},\\
&\|\nabla_\v\Delta_\v^{\frac{\alpha}{2}}\cT_{s,t}f(t)\|_\infty\leq C(t-s)^{-\frac{1}{\alpha}-1},\\
&\|\nabla_x\Delta_\v^{\frac{\alpha}{2}}\cT_{s,t}f(t)\|_\infty\leq C(t-s)^{-\frac{1}{\alpha}-2},
\end{align*}
we have for all $s\in[-1,1]$,
\begin{align*}
I_1&\leq C\int^\infty_2|\e^{\lambda(s-t)}-\e^{-\lambda t}|(t-s)^{-1}\dif t\\
&\leq C|\e^{\lambda s}-1|\int^\infty_2\e^{-\lambda t}\dif t=C|\e^{\lambda s}-1|\e^{-2\lambda}/\lambda\leq C,
\end{align*}
and for all $(s,x,\v)\in Q_1(0)$,
\begin{align*}
I_2(s,x,\v)\leq C\int^\infty_2\Big((t-s)^{-\frac{1}{\alpha}-1}+(t-s)^{-\frac{1}{\alpha}-2}\Big)\dif t\leq C.
\end{align*}
Moreover, by \eqref{NB6}, we have for all $s\in[0,1]$,
$$
I_3(s)\leq C.
$$
Combining the above calculations, we obtain \eqref{EP1} for $i=2$ with $C$ independent of $\lambda$.
For $i=1$, it is similar.
\end{proof}

Now we can give

\begin{proof}[Proof of Theorem \ref{Main} for $p\in(2,\infty)$]
By Lemmas \ref{Le34}, \ref{Le35} and \ref{Le37}, we know that
$$
\sP_i: L^\infty(\mR^{1+2d})\to BMO, i=1,2\mbox{ are bounded linear operators.}
$$
Estimate \eqref{BN77} for $p\in(2,\infty)$ follows by Theorem \ref{Th2} and the well-known estimate for $p=2$.
\end{proof}
\subsection{Proof of Theorem \ref{Main} for $p\in(1,2)$}

We shall use the dual argument to show that $\sP_i$, $i=1,2$ are still bounded linear operators in $L^p(\mR^{1+2d})$ for $p\in(1,2)$.
Let $\cT^*_{s,t}$ be the adjoint operator of $\cT_{s,t}$, that is,
$$
\int g\cT^*_{s,t} f=\int f\cT^*_{s,t} g.
$$
By definition \eqref{TST}, we have
$$
\cT^*_{s,t}f(x,\v):=\mE f\left(x+\int^t_s U_r\left[\v+\int^t_r\sigma_{r'}\dif L_{r',t}\right]\dif r,\v+\int^t_s\sigma_r\dif L_{r,t}\right).
$$
Let $p\in (1,2)$ and $q=\frac{p}{p-1}\in (2,\infty)$. By the dual relation between $L^p$ and $L^q$, we have
\begin{align*}
\|\sP_1f\|_p&=\sup_{h\in C^\infty_c(\mR^{1+2d)},\|h\|_q\leq 1}\int_{\mR^{1+2d}}\int^\infty_s\cT_{s,t}f\dif t\cdot \Delta^{\frac{\alpha}{2(1+\alpha)}}_x h\\
&=\sup_{h\in C^\infty_c(\mR^{1+2d)},\|h\|_q\leq 1}\int_{\mR^{1+2d}}f\cdot\left(\int^t_{-\infty}\cT^*_{s,t} \Delta^{\frac{\alpha}{2(1+\alpha)}}_x h\dif s\right).
\end{align*}
Since $\cT^*_{s,t} \Delta^{\frac{\alpha}{2(1+\alpha)}}_xh=\Delta^{\frac{\alpha}{2(1+\alpha)}}_x\cT^*_{s,t}h$, as in the previous subsection, one has
$$
\left\|\int^\cdot_{-\infty}\cT^*_{s,\cdot}\Delta^{\frac{\alpha}{2(1+\alpha)}}_x h\dif s\right\|_q=
\left\|\Delta^{\frac{\alpha}{2(1+\alpha)}}_x\int^\cdot_{-\infty}\cT^*_{s,\cdot} h\dif s\right\|_q\leq C\|h\|_q.
$$
Hence, by H\"older's inequality,
\begin{align}\label{UR2}
\|\sP_1f\|_p\leq C\|f\|_p.
\end{align}
Similarly, we have
\begin{align}\label{UR1}
\|\sP_2f\|_p=\sup_{h\in C^\infty_c(\mR^{1+2d)},\|h\|_q\leq 1}\int_{\mR^{1+2d}}f\cdot\left( \int^t_{-\infty}\cT^*_{s,t}\Delta^{\frac{\alpha}{2}}_\v h\dif s\right).
\end{align}
However, we can not treat it as $\sP_1$ because
$$
\cT^*_{s,t}\Delta^{\frac{\alpha}{2}}_\v h\not=\Delta^{\frac{\alpha}{2}}_\v\cT^*_{s,t} h.
$$
To overcome this difficulty, for $\eps\in(0,1)$, we introduce a new operator
\begin{align*}
\sQ_\eps f:&=\sQ^{\nu,\sigma, U}_\eps f(s,x,\v):=\int^t_{-\infty}\e^{\lambda(s-t)}\cT^{*,\nu,\sigma, U}_{s,t}\Delta_\v^{\frac{\alpha}{2}} f_\eps(s,x,\v)\dif s,
\end{align*}
where $f_\eps(t,x,\v)=f(t,\cdot)*\varrho_\eps(x,\v)$ so that $\sQ_\eps f$ is well defined for $f\in L^\infty(\mR^{1+2d})$.
Notice that $\sQ_0$ can be considered as the formal adjoint operator of $\sP_2$.
As in the previous subsection, we want to show that
$$
\sQ_\eps\mbox{ is a bounded linear operator from $L^\infty(\mR^{1+2d})$ to $BMO$.}
$$
First of all, as in Lemma \ref{Le34} we have
$$
\fint_{Q_r(t_0,x_0,\v_0)}\big|\sQ^{\nu,\sigma,U}_\eps f(s,x,\v)-a\big|^2=
\fint_{Q_1(0)}\big|\sQ^{\tilde\nu,\tilde\sigma,\tilde U}_\eps\tilde f(s,x,\v)-a\big|^2.
$$
where $\tilde\nu,\tilde\sigma,\tilde U$ and $\tilde f$ are defined as in Lemma \ref{Le34}.
We aim to prove that there is a constant $C=C(\kappa_0, p,d,\nu^{(\alpha)}_i,\alpha)>0$ independent of $\eps\in(0,1)$ such that for all $f\in L^\infty(\mR^{1+2d})$
with $\|f\|_\infty\leq 1$,
$$
\fint_{Q_1(0)}\big|\sQ^{\tilde\nu,\tilde\sigma,\tilde U}_\eps\tilde f(s,x,\v)-a\big|^2\leq C.
$$
Below we drop $\tilde\nu,\tilde\sigma,\tilde U$ and the tilde. As above, we make the following decomposition
$$
\sQ_\eps f=\left(\int^\cdot_{-2}+\int^{-2}_{-\infty}\right)\e^{\lambda(s-\cdot)}\cT^{*}_{s,\cdot}\Delta_\v^{\frac{\alpha}{2}} f_\eps(s)\dif s=:\sQ^\eps_{1}f+\sQ^\eps_{2}f.
$$
\bl
Let $\varphi\in C^\infty_c(\mR^{2d})$. For any $p\in[1,2]$, there exist constants $C_\varphi, \gamma>0$ such that
for all $h\in L^2(\mR^{2d})$ and $0<t-s\leq 3$,
\begin{align}\label{LK0}
\|\Delta^{\frac{\alpha}{2}}_\v(\cT_{s,t}(\varphi^2 h)-\varphi_{s,t}\cT_{s,t} (\varphi h))\|_p\leq C_\varphi(t-s)^{\gamma-1}\|h\|_2,
\end{align}
where $\varphi_{s,t}(x,\v):=\varphi(x+\Pi_{s,t}\v,\v)$ and $\Pi_{s,t}=\int^t_s U_r\dif r$.
\el
\begin{proof}
Let $p^\nu_{s,t}(x,\v)$ be the distributional density of $K^\nu_{s,t}$. Notice that
\begin{align*}
\nabla^2_\v\cT_{s,t}f(x,\v)
&=\int_{\mR^{2d}}f(x',\v')\nabla^2_\v p^{\nu}_{s,t}(x'-x-\Pi_{s,t}\cdot,\v'-\cdot)(\v)\dif x'\dif \v'\\
&=\int_{\mR^{2d}}f(x'+x+\Pi_{s,t}\v,\v'+\v)\Phi_{s,t}(x',\v')\dif x'\dif \v',
\end{align*}
where
$$
\Phi^{ij}_{s,t}=\sum_{i', j'}\Pi^{ii'}_{s,t}\Pi^{jj'}_{s,t}\p_{x_{i'}}\p_{x_{j'}}p^\nu_{s,t}+
2\sum_{i'}\Pi^{ii'}_{s,t}\p_{x_{i'}}\p_{\v_{j}}p^\nu_{s,t}+\p_{\v_{i}}\p_{\v_{j}}p^\nu_{s,t}.
$$
For any $\beta\in(0,\alpha)$, by \eqref{EV11}, it is easy to see that
\begin{align*}
&\|\nabla^2_\v\cT_{s,t}(\varphi^2 h)-\varphi_{s,t}\nabla^2_\v\cT_{s,t}(\varphi h)\|_p\\
&\quad\leq [\varphi]_\beta\|\varphi h\|_p\int_{\mR^{2d}}(|x'|^\beta+|\v'|^\beta)|\Phi_{s,t}(x',\v')|\dif x'\dif\v'\\
&\quad\leq C[\varphi]_\beta\|h\|_2(t-s)^{\frac{\beta(\alpha\wedge 1)-2}{\alpha}},
\end{align*}
where $[\varphi]_\beta:=\sup_{z\not=z'}|\varphi(z)-\varphi(z')|/|z-z'|^\beta$.
Furthermore, by the chain rule we have
\begin{align*}
\|\nabla^2_\v(\cT_{s,t}(\varphi^2 h)-\varphi_{s,t}\cT_{s,t} (\varphi h))\|_p
\leq C_\varphi\|h\|_2(t-s)^{\frac{\beta(\alpha\wedge 1)-2}{\alpha}}.
\end{align*}
Hence, by definition \eqref{Def3} and \eqref{EV31}, we have
\begin{align*}
&\|\Delta^{\frac{\alpha}{2}}_\v(\cT_{s,t}(\varphi^2 h)-\varphi_{s,t}\cT_{s,t}(\varphi h))\|_p\\
&\qquad\leq\int_{\mR^d}\|\delta^{(2)}_\v(\cT_{s,t}(\varphi^2 h)-\varphi_{s,t}\cT_{s,t} (\varphi h))\|_p\frac{\dif\v}{|\v|^{d+\alpha}}\\
&\qquad\leq 4\int_{|\v|>(t-s)^{(2-\beta(\alpha\wedge 1))/(2\alpha)}}\Big(\|\cT_{s,t}(\varphi^2 h)\|_p+\|\varphi_{s,t}\cT_{s,t} (\varphi h)\|_p\Big)\frac{\dif\v}{|\v|^{d+\alpha}}\\
&\qquad\quad+C_\phi\|h\|_2(t-s)^{\frac{\beta(\alpha\wedge 1)-2}{\alpha}}\int_{|\v|\leq (t-s)^{(2-\beta(\alpha\wedge 1))/(2\alpha)}}|\v|^2\frac{\dif\v}{|\v|^{d+\alpha}}\\
&\qquad\leq C_\phi\|h\|_2(t-s)^{\frac{\beta(\alpha\wedge 1)}{2}-1}.
\end{align*}
Thus, we obtain \eqref{LK0}.
\end{proof}

\bl
Under  \eqref{Kapp} and \eqref{Con2}, there is a positive constant $C$ only depending on $\kappa_0, p,d,\nu^{(\alpha)}_i,\alpha$
such that for all $f\in L^\infty(\mR^{1+2d})$ with $\|f\|_\infty\leq 1$ and all $\eps\in(0,1)$,
\begin{align}
\int_{Q_1(0)}|\sQ^*_{\eps} f(s,x,\v)|^2\leq C.\label{BN4'}
\end{align}
\el
\begin{proof}
For $t\in\mR$, define
$$
u(t,x,\v):=\int^t_{-\infty}\e^{\lambda(s-t)}\cT^*_{s,t}\Delta^{\frac{\alpha}{2}}_\v ((1_{[-2,1]}f_\eps)(s))(x,\v)\dif s.
$$
Let $\varphi$ be a nonnegative smooth cutoff function in $\mR^{2d}$ with $\varphi(x,\v)=1$ for $|(x,\v)|\leq 4$ and  $\varphi(x)=0$ for $|(x,\v)|>8$. We have
\begin{align*}
\|u\|_{L^2(Q_1(0))}&\leq\|u\varphi^2\|_{2}=\sup_{h\in C^\infty_c(\mR^{1+2d}), \|h\|_2\leq 1}\int_{\mR^{1+2d}}u\varphi^2h\\
&=\sup_{h\in C^\infty_c(\mR^{1+2d}), \|h\|_2\leq 1}\int_{\mR^{1+2d}}1_{[-2,1]}(s)f_\eps\Delta^{\frac{\alpha}{2}}_\v
\int^\infty_s \e^{\lambda(s-t)}\cT_{s,t}(\varphi^2 h (t))\dif t\\
&\leq\|f\|_\infty\sup_{h\in C^\infty_c(\mR^{1+2d}), \|h\|_2\leq 1}
\left\|1_{[-2,1]}\Delta^{\frac{\alpha}{2}}_\v\int^1_{\cdot}\e^{\lambda(s-t)}\cT_{s,t}(\varphi^2 h (t))\dif t\right\|_1.
\end{align*}
Let $h\in C^\infty_c(\mR^{1+2d})$ with $\|h\|_2\leq 1$. By \eqref{LK0}, we have
\begin{align*}
&\left\|1_{[-2,1]}\Delta^{\frac{\alpha}{2}}_\v\int^1_{\cdot}\e^{\lambda(s-t)}\cT_{s,t}(\varphi^2 h(t))\dif t\right\|_1\\
&\quad\leq\left\|1_{[-2,1]}\Delta^{\frac{\alpha}{2}}_\v\int^1_{\cdot}\e^{\lambda(s-t)}\varphi_{s,t}\cT_{s,t}(\varphi h(t))\dif t\right\|_1+C.
\end{align*}
Since $\varphi_{s,t}=\varphi(x+\Pi_{s,t}\v,\v)$ has support $\Big\{(x,\v): |(x,\v)|\leq 8(\|U\|_\infty+1)\Big\}$, and for any $\gamma\in(0,\alpha\wedge 1)$,
\begin{align*}
\|\delta^{(1)}_\v\cT_{s,t}(\varphi h)\|_2&\leq (\|\nabla_\v\cT_{s,t}(\varphi h)\|_2|\v|)\wedge (2\|\cT_{s,t}(\varphi h)\|_2)\\
&\leq C((t-s)^{-\frac{1}{\alpha}}|\v|)\wedge 1\leq C(t-s)^{-\frac{\gamma}{\alpha}}|\v|^\gamma,
\end{align*}
by \eqref{EV441} and \eqref{LK0}, we have
\begin{align*}
&\Bigg\|1_{[-2,1]}\Delta^{\frac{\alpha}{2}}_\v\int^1_{\cdot}\e^{\lambda(s-t)}\varphi_{s,t}\cT_{s,t}(\varphi h(t))\dif t\Bigg\|_1\\
&\quad\leq\left\|1_{[-2,1]}\int^1_{\cdot}\e^{\lambda(s-t)}\varphi_{s,t}\Delta^{\frac{\alpha}{2}}_\v\cT_{s,t}(\varphi h(t))\dif t\right\|_1+C\\
&\quad\leq\left\|1_{[-2,1]}\int^1_{\cdot}\e^{\lambda(s-t)}\varphi_{s,t}\Delta^{\frac{\alpha}{2}}_\v\cT_{s,t}(\varphi h(t))\dif t\right\|_2+C\\
&\quad\leq\left\|1_{[-2,1]}\Delta^{\frac{\alpha}{2}}_\v\int^1_{\cdot}\e^{\lambda(s-t)}\varphi_{s,t}\cT_{s,t}(\varphi h(t))\dif t\right\|_2+C\\
&\quad\leq\left\|1_{[-2,1]}\Delta^{\frac{\alpha}{2}}_\v\int^1_{\cdot}\e^{\lambda(s-t)}\cT_{s,t}(\varphi^2 h (t))\dif t\right\|_2+C\\
&\quad\leq C\|\varphi^2 h\|_2+C\leq C.
\end{align*}
Combining the above calculations, we obtain \eqref{BN4'}.
\end{proof}
The following lemma is the same as in Lemma \ref{Le37}.
\bl
Under  \eqref{Kapp} and \eqref{Con2}, there is a positive	constant $C$ only depending on $\kappa_0, p,d,\nu^{(\alpha)}_i,\alpha$
such that for all $f\in L^\infty(\mR^{1+2d})$ with $\|f\|_\infty\leq 1$ and $\eps\in(0,1)$,
$$
\int_{Q_1(0)}|\sQ^\eps_{2}f(t,x,\v)-\sQ^\eps_{2}f(0,0,0)|^2\leq C.
$$
\el
\begin{proof}
By \eqref{JH2} with $\beta=\gamma=\alpha$, we have for all $t\in[-1,1]$,
\begin{align*}
&\int^{-2}_{-\infty}\big|\cT^*_{s,t}\Delta_\v^{\frac{\alpha}{2}}f_\eps(s,0,0)
-\cT^*_{s,0}\Delta_\v^{\frac{\alpha}{2}}f_\eps(s,0,0)\big|\dif s\\
&\quad\leq\int^{-2}_{-\infty}\!\int^t_0\big|\p_r\cT^*_{s,r}\Delta_\v^{\frac{\alpha}{2}}f_\eps(s,0,0)\big|\dif r\dif s\\
&\quad=\int^{-2}_{-\infty}\!\int^t_0\big|\sL^{\nu_r}_{\sigma_r, \v}\cT^*_{s,r}\Delta_\v^{\frac{\alpha}{2}}f_\eps(s,0,0)\big|\dif r\dif s\\
&\quad\leq C\int^{-2}_{-\infty}\!\int^t_0(t-r)^{-2}\dif r\dif t\leq C.
\end{align*}
Using this estimate, as in the proof of Lemma \ref{Le37}, we obtain the desired estimate.
\end{proof}
Now we can give

\begin{proof}[Proof of Theorem \ref{Main} for $p\in(1,2)$]
By Lemmas \ref{Le34}, \ref{Le35} and \ref{Le37}, we know that
$$
\sQ_\eps: L^\infty(\mR^{1+2d})\to BMO\mbox{ is bounded with norm independent of $\eps$.}
$$
Moreover, by duality, we also have
$$
\sQ_\eps: L^2(\mR^{1+2d})\to L^2(\mR^{1+2d})\mbox{ is bounded with norm independent of $\eps$.}
$$
Hence, for $q=p/(p-1)\in(2,\infty)$, by Theorem \ref{Th2}, we have for some $C>0$ independent of $\eps$,
$$
\|\sQ_\eps f\|_q=\left\|\int^t_{-\infty}\e^{\lambda(s-t)}\cT^*_{s,t}\Delta^{\frac{\alpha}{2}}_\v f_\eps\dif s\right\|_q\leq C\|f\|_q.
$$
Now going back to \eqref{UR1}, for $p\in(1,2)$, by Fatou's lemma, we get
\begin{align*}
\|\sP_2f\|_p&\leq \|f\|_p\sup_{h\in C^\infty_c(\mR^{1+2d)},\|h\|_q\leq 1}\left\|\int^t_{-\infty}\e^{\lambda(s-t)}\cT^*_{s,t}\Delta^{\frac{\alpha}{2}}_\v h\dif s\right\|_q\\
&\leq \|f\|_p\sup_{h\in C^\infty_c(\mR^{1+2d)},\|h\|_q\leq 1}\varliminf_{\eps\to 0}
\left\|\int^t_{-\infty}\e^{\lambda(s-t)}\cT^*_{s,t}\Delta^{\frac{\alpha}{2}}_\v h_\eps\dif s\right\|_q\leq C\|f\|_p,
\end{align*}
which together with \eqref{UR2} gives \eqref{BN77} for $p\in(1,2)$.
\end{proof}

\section{Appendix}
\subsection{Carleman's representation for Boltzmann's equation}
Let us first show the following elementary formula in calculus.
\bl
We have
\begin{align}\label{For}
\int_{\mR^d}\!\!\int_{\mS^{d-1}}F(x,\omega)\dif \omega\dif x
=
\int_{\mR^d}\!\!\int_{\{h\cdot w=0\}}F(h\pm w,\bar w)|w|^{1-d}\dif h\dif w,
\end{align}
where $\bar w:=w/|w|$ and we have used the convention that $F(h\pm w,\bar w)=F(h+w,\bar w)+F(h-w,\bar w)$. In particular,
if $F(x,\omega)=F(x,-\omega)$, then
$$
\int_{\mR^d}\!\!\int_{\mS^{d-1}}F(x,\omega)\dif \omega\dif x=2\int_{\mR^d}\!\!\int_{\{h\cdot w=0\}}F(h+w,\bar w)|w|^{1-d}\dif h\dif w.
$$
\el
\begin{proof}
By the co-area formula and the change of variables, we have
\begin{align*}
\int_{\mS^{d-1}}\!\!\int_{\mR^d}F(x,\omega)\dif x\dif \omega
&=\int_{\mS^{d-1}}\!\!\int^\infty_0\!\!\int_{\{\<h,\omega\>=\pm r\}}F(h,\omega)\dif h\dif r\dif \omega\\
&=\int_{\mS^{d-1}}\!\!\int^\infty_0\!\!\int_{\{\<h\pm r\omega,\omega\>=0\}}F(h,\omega)\dif h\dif r\dif \omega\\
&=\int_{\mS^{d-1}}\!\!\int^\infty_0\!\!\int_{\{\<h,\omega\>=0\}}F(h\pm r\omega,\omega)\dif h\dif r\dif \omega\\
&=\int_{\mR^{d}}\!\!\int_{\{\<h,w\>=0\}}F(h\pm w,\bar w)|w|^{1-d}\dif h\dif w.
\end{align*}
The desired formula follows.
\end{proof}

By a change of variables and
\eqref{For},
noting that $\<h,w\>=0$,  one can rewrite the collision operator $Q(f,g)$ as:
\begin{align*}
Q(f,g)(\v)&=\int_{\mR^d}\int_{\mS^{d-1}}\Big[f(\v-\v_*+\<\v_*,\omega\>\omega)g(\v-\<\v_*,\omega\>\omega)\\
&\qquad\qquad\qquad-f(\v-\v_*)g(\v)\Big]B(|\v_*|,\omega)\dif\omega\dif \v_*\\
&=2\int_{\mR^d}\int_{\{h\cdot w=0\}}\Big[f(\v-h-w+\<h+w,\bar w\>\bar w)g(\v-\<h+w,\bar w\>\bar w)\\
&\qquad\qquad\qquad-f(\v-h- w)g(\v)\Big]B(|h+ w|,\bar w)|w|^{1-d}\dif h\dif w\\
&=2\int_{\mR^d}\int_{\{h\cdot w=0\}}\Big[f(\v-h)g(\v-w)-f(\v-h-w)g(\v)\Big]\\
&\qquad\qquad\qquad\qquad\times B(|h+w|,\bar w)|w|^{1-d}\dif h\dif w,
\end{align*}
which gives representation \eqref{Bol} by changing $w$ into $-w$.
\subsection{Proof of Theorem \ref{Th2}}
Let us introduce a quasi-metric in $\mR^{1+2d}$ as follows:
\begin{align*}
&\rho((t_0,x_0,\v_0),(t_1,x_1,\v_1))\\
&:=|t_0-t_1|^{\frac{1}{\alpha}}+|\v_0-\v_1|+|x_0-x_1+\Pi_{t_0,t_1}\v_1|^{\frac{1}{1+\alpha}}+|x_1-x_0+\Pi_{t_1,t_0}\v_0|^{\frac{1}{1+\alpha}},
\end{align*}
where $\Pi_{t_0, t_1}:=\int^{t_1}_{t_0}U_r\dif r$. More precisely, $\rho$ satisfies
\begin{enumerate}[(i)]
\item $\rho((t_0,x_0,\v_0),(t_1,x_1,\v_1))=0\Rightarrow t_0=t_1, x_0=x_1, \v_0=\v_1$.
\item $\rho((t_0,x_0,\v_0),(t_1,x_1,\v_1))=\rho((t_1,x_1,\v_1),(t_0,x_0,\v_0))$.
\item For some constant $c_0\geq 1$ and any points $(t_i,x_i,\v_i)\in\mR^{1+2d}, i=0,1,2$, it holds that
\begin{align*}
&\rho((t_0,x_0,\v_0),(t_2,x_2,\v_2))\\
&\quad\leq c_0\Big(\rho((t_0,x_0,\v_0),(t_1,x_1,\v_1))+\rho((t_1,x_1,\v_1),(t_2,x_2,\v_2))\Big).
\end{align*}
\end{enumerate}
Given $(t_0,x_0,\v_0)\in\mR^{1+2d}$ and $r>0$, a ``ball'' in $\mR^{1+2d}$ with radius $r$ with respect to the quasi-metric $\rho$ is defined by
$$
\widetilde Q_r(t_0,x_0,\v_0):=\Big\{(t,x,\v)\in\mR^{1+2d}: \rho((t_0,x_0,\v_0), (t,x,\v))<r\Big\}.
$$
Recalling the definition of  the ``ball'' $Q_r$ in (\ref{Ball}),
we have the following relation between  $\widetilde Q_r$ and $Q_r$,
whose proof is obvious by definitions.
\bl\label{Le51}
Let $c_1:=(4+\|U\|_\infty)^\alpha$. For any $r>0$ and $(t_0,x_0,\v_0)\in\mR^{1+2d}$, we have
\begin{align}
\widetilde Q_r(t_0,x_0,\v_0)\subset Q_r(t_0,x_0,\v_0)\subset\widetilde Q_{c_1r}(t_0,x_0,\v_0).
\end{align}
In particular, let  $c_2:=(2c_1)^{1+(2+\alpha)d}$, the following doubling property holds:
\begin{align}
|\widetilde Q_{2r}(t_0,x_0,\v_0)|\leq c_2|\widetilde Q_{r}(t_0,x_0,\v_0)|.\label{EY1}
\end{align}
\el

The doubling property (\ref{EY1}) means that $(\mR^{1+2d},\rho,\dif x)$ is a space of homogenous type in the sense of \cite[Definition 1]{Ch}.
Thus by the $T(b)$ theorem (see \cite[Theorem 11]{Ch}) ,
we have \bl\label{Le42}
With respect to the space $(\mR^{1+2d},\rho,\dif x)$,
there exists a collection of open subsets $\{O_{nj}\subset\mR^{1+2d}, n\in\mZ, j\in \cI_n\}$, where $\cI_n$ denotes some index set depending on $n$,
and constants $\delta\in(0,1)$, $a_0>0$ and $c_3>0$ such that
\begin{enumerate}[(i)]
\item For each $n\in\mN$, $|\mR^{1+2d}\setminus(\cup_{j\in\mI_n})O_{nj})|=0$.
\item If $n\geq k$, then either $O_{nj}\subset O_{ki}$ or $O_{nj}\cap O_{ki}=\emptyset$.
\item For each $(n,j)$ and $k<n$, there is a unique $i\in \cI_k$ such that $O_{nj}\subset O_{ki}$.
\item Diameter of $O_{kj}$ is less than $c_3\delta^n$, and hence, for each $(t,x,\v)\in O_{nj}$, we have
$$
O_{nj}\subset \widetilde Q_{\delta^n}(t,x,\v).
$$
\item For each $(n,j)$, $O_{nj}$ contains some ball $\widetilde Q_{a_0\delta^n}$, and so
$$
|O_{nj}|\geq c_3\delta^n.
$$
\end{enumerate}
\el
Let $O_{nj}$ be as in the above lemma, which will play the role of ``cube'' in the classical Calder\'on-Zygmund's decomposition. We write
$$
\mC_n:=\Big\{O_{nj}, j\in \cI_n\Big\}.
$$
If we define $O:=\cap_{n\in\mZ}(\cup_{j\in\mI_n})O_{nj})$, then by (i), the complement $O^c$ has null Lebesgue measure.
By restricting on $O$, without loss of generality,
we may assume that $\mC_n$ is a partition of $\mR^{1+2d}$. Thus, for each $(t,x,\v)\in\mR^{1+2d}$, there is a unique $O_n\in\mC_n$ such that
$(t,x,\v)\in O_n$. We will also denote this $O_n$ by $O_n(t,x,\v)$, and for any local integrable function $f$, define
$$
f_{|_n}(t,x,\v):=\fint_{O_n(t,x,\v)}f(t',x',\v')\dif \v'\dif x'\dif t'.
$$
The function $f_{|_n}$ can be considered as a ``conditional function of $f$ given $\mC_n$''.
By Lemma \ref{Le42}, one sees that $\{\mC_n,n\in\mN\}$ forms a sequence of partitions in the sense of \cite[Definition 1, p.74]{Kr}. More precisely,
\begin{enumerate}[(i)]
\item For each $n$ and $O_n\in\mC_n$, there is a unique $O_{n-1}\in\mC_{n-1}$ such that $O_n\subset O_{n-1}$, and
$$
|O_{n-1}|\leq C_{d,\alpha}|O_n|.
$$
\item For any continuous function $f$ on $\mR^{1+2d}$, we have
$$
\lim_{n\to\infty}f_{|_n}(t,x,\v)\to f(t,x,\v),\ \forall (t,x,\v)\in\mR^{1+2d}.
$$
\end{enumerate}
Now we can give
\begin{proof}[Proof of Theorem \ref{Th26}]
We define another sharp function associated to $\{\mC_n, n\in\mZ\}$ by
$$
\tilde\cM^\sharp f(t,x,\v):=\max_{n\in\mZ}\fint_{O_n(t,x,\v)}|f-f_{|_n}|.
$$
By \cite[Theorem 10, p.81]{Kr}, for any $p\in(1,\infty)$, we have
\begin{align}
\|f\|_p\leq C\|\tilde\cM^\sharp f\|_p.\label{EN2}
\end{align}
On the other hand, by (iv) and (v) of Lemma \ref{Le42}, we have
$$
O_n(t,x,\v)\subset \widetilde Q_{c\delta^n}(t,x,\v)\subset Q_{c\delta^n}(t,x,\v)
$$
and
$$
|O_n(t,x,\v)|\geq c\delta^n\geq |Q_{c\delta^n}(t,x,\v)|.
$$
Therefore,
\begin{align}
\fint_{O_n(t,x,\v)}|f-f_{|_n}|&\leq\fint_{O_n(t,x,\v)}\fint_{O_n(t,x,\v)}|f(t',x',\v')-f(t'',x'',\v'')|\no\\
&\leq\fint_{Q_{c\delta^n}(t,x,\v)}\fint_{Q_{c\delta^n}(t,x,\v)}|f(t',x',\v')-f(t'',x'',\v'')|\no\\
&\leq 2\fint_{Q_{c\delta^n}(t,x,\v)}|f-f_{Q_{c\delta^n}(t,x,\v)}|\leq\cM^\sharp f(t,x,\v).\label{EN3}
\end{align}
Estimate (\ref{EU4}) now follows by (\ref{EN2}) and (\ref{EN3}).
\end{proof}

\subsection{Proof of Theorem \ref{Main} for $p=2$}

In this subsection we give a proof of Theorem \ref{Main} for $p=2$. Let us first recall a key estimate due to Bouchut \cite{Bo}.
Since there is a time inhomogeneous matrix $U_s$ in our formulation, we need to modify the proof given in \cite{Al}.

\bt\label{Th44}
Let $U:\mR\to\mM^d$ satisfy
$$
\kappa_0:=\sup_{s}\|U_s\|+\sup_{s<t}\Big((t-s)\|\Pi_{s,t}^{-1}\|\Big)<\infty,
$$
where $\Pi_{s,t}:=\int^t_sU_r\dif r$.
Let $u,f\in L^2(\mR^{1+2d})$ with
$\Delta_\v^{\frac{\alpha}{2}} u\in L^2(\mR^{1+2d})$
for some $\alpha\geq 0$, and satisfy
\begin{align}
\p_s u+U_s\v\cdot\nabla_x u+f=0
\quad \hbox{in the distributional sense.}
\label{EW1}
\end{align}
Then for some $C=C(d,\alpha,\kappa_0)>0$, we have
\begin{align}
\|\Delta^{\frac{\alpha}{2(1+\alpha)}}_xu\|_2\leq C\|\Delta^{\frac{\alpha}{2}}_\v u\|^{\frac{1}{1+\alpha}}_2\|f\|^{\frac{\alpha}{1+\alpha}}_2.\label{EV10}
\end{align}
\et

\begin{proof}
We follow the argument of \cite{Al}
with modification to deal with the time-dependent case.
Taking Fourier transform in $(x, \v)$-variables on both sides of \eqref{EW1},
we have
\begin{align}\label{EM3}
\p_s \hat u-U^*_s\xi\cdot\nabla_\eta \hat u+\hat f=0 .
\end{align}
Here
$$
\hat u(t, \xi,\eta)=\int_{\mR^{2d}}\e^{{\rm i}(\xi\cdot x+\eta\cdot\v)}u(t, x,\v)\dif x\dif\v,
$$
and $\hat f(t, \xi, \eta)$ is defined in a similar way.
Let
$\psi: [0,\infty)\to[0,1]$ be a smooth function with $\psi (s)=1$ for $s<1$ and $\psi (s)=0$ for $s>2$. For $\eps\in(0,1)$, define
$$
\phi_\eps(\xi,\eta):=\psi (\eps|\eta|/|\xi|^{1/(1+\alpha)}).
$$
By Planchel's identity, we have
\begin{align}
\begin{split}
\|\Delta^{\frac{\alpha}{2(1+\alpha)}}_xu\|^2_2
\, =& \int^\infty_{-\infty}\!\int_{\mR^{2d}}|\xi|^{\frac{2\alpha}{1+\alpha}}|\hat u|^2(s,\xi,\eta)\dif\eta\dif\xi\dif s\\
\leq & 2\int^\infty_{-\infty}\!\int_{\mR^{2d}}|\xi|^{\frac{2\alpha}{1+\alpha}}(1-\phi_\eps(\xi,\eta))^2|\hat u|^2(s,\xi,\eta)\dif\eta\dif\xi\dif s\\
&+2\int^\infty_{-\infty}\!\int_{\mR^{2d}}|\xi|^{\frac{2\alpha}{1+\alpha}}\phi^2_\eps(\xi,\eta)|\hat u|^2(s,\xi,\eta)\dif\eta\dif\xi\dif s \\
=:& \sI_\eps+\sJ_\eps.
\end{split}
\label{EM1}
\end{align}
For $\sI_\eps$, by the definition of $\phi_\eps$, we have
\begin{align}\label{EM2}
\sI_\eps\leq 2\eps^{2\alpha}\int^\infty_{-\infty}\!\int_{\mR^{2d}}|\eta|^{2\alpha}|\hat u|^2(s,\xi,\eta)\dif\eta\dif\xi\dif s=2\eps^{2\alpha}\|\Delta^{\frac{\alpha}{2}}_\v u\|^2_2.
\end{align}
To treat $\sJ_\eps$, let us write
$$
\hat u_{\eps}:=\phi_\eps\hat u,\ \ \hat f_\eps:=\phi_\eps\hat f,\ \ g_\eps:=\hat f_\eps+(U^*_s\xi\cdot\nabla_\eta\phi_\eps) \hat u.
$$
Then by \eqref{EM3}, it is easy to see that
$$
\p_s \hat u_\eps-U^*_s\xi\cdot\nabla_\eta \hat u_\eps+\hat g_\eps=0.
$$
Multiplying both sides by the complex conjugate of $\hat u_\eps$, we obtain
$$
\p_s |\hat u_\eps|^2-U^*_s\xi\cdot\nabla_\eta |\hat u_\eps|^2+2\Re\mathrm{e}(\hat g_\eps,\bar{\hat u}_\eps)=0.
$$
It follows that
\begin{align*}
|\hat u_\eps|^2(s,\xi,\eta)
&= - \int_s^\infty \frac{\dif }{\dif t}  |\hat u_\eps|^2(t,\xi, \eta -\Pi^*_{s, t} \xi) \dif t\\
&= 2\int^\infty_s\Re\mathrm{e}(\hat g_\eps,\bar{\hat u}_\eps)(t, \xi,\eta-\Pi^*_{s,t}\xi)\dif t.
\end{align*}
Since the support of $\phi_\eps$ is contained in $\big\{\eps|\eta|<2|\xi|^{1/(1+\alpha)}\big\}$, we get
\begin{align*}
\sI_\eps&=\int^\infty_{-\infty}\!\int_{\mR^{2d}}|\xi|^{\frac{2\alpha}{1+\alpha}}1_{\{\eps|\eta|<2|\xi|^{1/(1+\alpha)}\}}|\hat u_\eps|^2(s,\xi,\eta)\dif\eta\dif\xi\dif s\\
&\leq2\int^\infty_{-\infty}\!\int_{\mR^{2d}}|\xi|^{\frac{2\alpha}{1+\alpha}}\int^\infty_s1_{\{\eps|\eta|<2|\xi|^{1/(1+\alpha)}\}}
|\hat g_\eps|\,|\hat u_\eps|(t, \xi,\eta-\Pi^*_{s,t}\xi)\dif t\dif\eta\dif\xi\dif s\\
&=2\int^\infty_{-\infty}\!\int_{\mR^{2d}}|\xi|^{\frac{2\alpha}{1+\alpha}}\int^\infty_s1_{\{\eps|\eta-\Pi^*_{s,t}\xi|<2|\xi|^{1/(1+\alpha)}\}}
|\hat g_\eps|\,|\hat u_\eps|(t, \xi,\eta)\dif t\dif\eta\dif\xi\dif s\\
&=2\int^\infty_{-\infty}\!\int_{\mR^{2d}}|\xi|^{\frac{2\alpha}{1+\alpha}}\left(\int^t_{-\infty}1_{\{\eps|\eta-\Pi^*_{s,t}\xi|<2|\xi|^{1/(1+\alpha)}\}}
\dif s\right)|\hat g_\eps|\,|\hat u_\eps|(t, \xi,\eta)\dif\eta\dif\xi\dif t.
\end{align*}
Let us estimate the integral in the bracket. By the assumption, we have
\begin{align*}
&\int^t_{-\infty}1_{\{\eps|\eta-\Pi^*_{s,t}\xi|<2|\xi|^{1/(1+\alpha)}\}}\dif s
\leq\int^t_{-\infty}1_{\{|\Pi^*_{s,t}\xi|\leq |\eta|+2|\xi|^{1/(1+\alpha)}/\eps\}}\dif s\\
&\quad\leq\int^t_{-\infty}1_{\{\kappa^{-1}_0(t-s)|\xi|\leq |\eta|+2|\xi|^{1/(1+\alpha)}/\eps\}}\dif s
=\kappa_0\Big(|\eta|+2|\xi|^{1/(1+\alpha)}/\eps\Big)/|\xi|.
\end{align*}
Moreover, by the definition of $\hat g_\eps$, we also have
\begin{align*}
|\hat g_\eps||\hat u_\eps|&\leq \left(|\hat f_\eps|+\kappa_0\eps\|\psi'\|_\infty|\xi|^{\alpha/(1+\alpha)}
1_{\{1\leq\eps |\eta|/|\xi|^{1/(1+\alpha)}\leq 2\}}|\hat u|\right)|\hat u_\eps|.
\end{align*}
Therefore,
\begin{align*}
\sJ_\eps&\leq \frac{2\kappa_0}{\eps}\int^\infty_{-\infty}\!\int_{\mR^{2d}}
|\xi|^{\frac{\alpha-1}{1+\alpha}}(\eps|\eta|+2|\xi|^{1/(1+\alpha)})|\hat f_\eps|
|\hat u_\eps|(t, \xi,\eta)\dif\eta\dif\xi\dif t\\
&\quad+2\kappa^2_0\|\psi'\|_\infty\int^\infty_{-\infty}\!\int_{\mR^{2d}}
|\xi|^{\frac{2\alpha-1}{1+\alpha}}(\eps|\eta|+2|\xi|^{1/(1+\alpha)})\\
&\qquad\times 1_{\{1\leq\eps |\eta|/|\xi|^{1/(1+\alpha)}\leq 2\}}|\hat u|\,|\hat u_\eps|(t, \xi,\eta)\dif\eta\dif\xi\dif t\\
&\leq \frac{8\kappa_0}{\eps}\int^\infty_{-\infty}\!\int_{\mR^{2d}}
|\xi|^{\frac{\alpha}{1+\alpha}}|\hat f|\,|\hat u|(t, \xi,\eta)\dif\eta\dif\xi\dif t\\
&\quad+8\kappa^2_0\|\psi'\|_\infty\int^\infty_{-\infty}\!\int_{\mR^{2d}}
(\eps|\eta|)^{2\alpha}|\hat u|^2(t, \xi,\eta)\dif\eta\dif\xi\dif t\\
&\leq C\eps^{-2}\|f\|^2_2+\tfrac{1}{2}\|\Delta^{\frac{\alpha}{2(1+\alpha)}}_xu\|^2_2+8\kappa^2_0\|\psi'\|_\infty\eps^{2\alpha}\|\Delta^{\frac{\alpha}{2}}_\v u\|^2_2.
\end{align*}
 Combining this with \eqref{EM1} and \eqref{EM2}, we obtain
\begin{align*}
\|\Delta^{\frac{\alpha}{2(1+\alpha)}}_xu\|^2_2&\leq(1+8\kappa^2_0\|\psi'\|_\infty)\eps^{2\alpha}\|\Delta^{\frac{\alpha}{2}}_\v u\|^2_2
+C\eps^{-2}\|f\|^2_2+\tfrac{1}{2}\|\Delta^{\frac{\alpha}{2(1+\alpha)}}_xu\|_2^2,
\end{align*}
which gives the desired estimate by letting $\eps=(C\|f\|^2_2/\|\Delta^{\frac{\alpha}{2}}_\v u\|^2_2)^{\frac{1}{2(1+\alpha)}}$.
\end{proof}
Now we can give

\begin{proof}[Proof of Theorem \ref{Main} for $p=2$]

Without loss of generality, we may assume $f\in C^\infty_c(\mR^{1+2d})$.
It follows from   Fourier transformation, \eqref{TST0}, \eqref{NB9}, and H\"older's inequality that
\begin{align*}
&\| \Delta^{\alpha/2}_v u^\lambda \|_2
= \int^\infty_{-\infty}\left\|\Delta^{{\alpha}/{2}}_\v\int^\infty_s\e^{\lambda(s-t)}\cT_{s,t}f(t,\cdot,\cdot)\dif t\right\|_2^2\dif s\\
&=\int^\infty_{-\infty}\int_{\mR^{2d}}|\eta|^{2\alpha}\left|\int^\infty_s\e^{\lambda(s-t)}\widehat{\cT_{s,t}f}(t,\xi,\eta)\dif t\right|^2\dif\eta\dif\xi\dif s\\
&=\int^\infty_{-\infty}\int_{\mR^{2d}}|\eta|^{2\alpha}\left|\int^\infty_s\e^{\lambda(s-t)}
\e^{-\int^t_s\psi^{\nu_r}_{\sigma_r}(\Pi^*_{s,r}\xi-\eta)\dif r}\hat f(t,\xi,\eta-\Pi^*_{s,r}\xi )
\dif t\right|^2\dif\eta\dif\xi\dif s\\
&\leq\int^\infty_{-\infty}\int_{\mR^{2d}}\left(\int^\infty_s|\eta|^{\alpha}
\e^{-\int^t_s\psi^{\nu_r}_{\sigma_r}(\Pi^*_{s,r}\xi-\eta)\dif r}|\hat f(t,\xi,\eta-\Pi^*_{s,t}\xi)|^2\dif t\right)\\
&\qquad\qquad\times\left(\int^\infty_s|\eta|^{\alpha}
\e^{-\int^t_s\psi^{\nu_r}_{\sigma_r}(\Pi^*_{s,r}\xi-\eta)\dif r}\dif t\right)\dif\eta\dif\xi\dif s.
\end{align*}
By \eqref{Con1} and a similar argument as that for \eqref{NB2}, we have
$$
\int^t_s\psi^{\nu_r}_{\sigma_r}(\Pi^*_{s,r}\xi-\eta)\dif r\geq\kappa_1\int^t_s|\Pi^*_{s,r}\xi-\eta|^\alpha\dif r
\geq c_1(t-s)|(\eta,(t-s)\xi)|^\alpha,
$$
and
$$
\int^t_s\psi^{\nu_r}_{\sigma_r}(\Pi^*_{r,t}\xi+\eta)\dif r\geq\kappa_1\int^t_s|\Pi^*_{r,t}\xi+\eta|^\alpha\dif r
\geq c_1(t-s)|(\eta,(t-s)\xi)|^\alpha.
$$
Hence,
$$
\int^\infty_s|\eta|^{\alpha}\e^{-\int^t_s\psi^{\nu_r}_{\sigma_r}( \Pi^*_{s,r}\xi-\eta)\dif r}\dif t\leq
\int^\infty_s|\eta|^{\alpha}
\e^{-c_1(t-s)|\eta|^\alpha}\dif t=\frac{1}{c_1},
$$
and
\begin{align*}
&\int^t_{-\infty}|\eta+\Pi^*_{s,t}\xi|^{\alpha}\e^{-\int^t_s\psi^{\nu_r}_{\sigma_r}(\Pi^*_{r,t}\xi+\eta)\dif r}\dif s
\leq 2\int^t_{-\infty}|\eta|^{\alpha}\e^{-c_1(t-s)|\eta|^\alpha}\dif s\\
&\qquad+2\|U\|_\infty^\alpha\int^t_{-\infty}((t-s)|\xi|)^{\alpha}\e^{-c_1(t-s)^{1+\alpha}|\xi|^\alpha}\dif s=\frac{2}{c_1}+\frac{2\|U\|^\alpha_\infty}{c_1(1+\alpha)}.
\end{align*}
Thus, by the change of variables and Fubini's theorem, we further have
\begin{align*}
&\int^\infty_{-\infty}\left\|\Delta^{\frac{\alpha}{2}}_\v\int^\infty_s\e^{\lambda(s-t)}\cT_{s,t}f(t,\cdot,\cdot)\dif t\right\|_2^2\dif s\\
&\leq\frac{1}{c_1}\int^\infty_{-\infty}\int_{\mR^{2d}}\left(\int^\infty_s|\eta|^{\alpha}
\e^{-\int^t_s\psi^{\nu_r}_{\sigma_r}(\Pi^*_{s,r}\xi-\eta)\dif r}|\hat f(t,\xi,\eta-\Pi^*_{s,t}\xi)|^2\dif t\right)\dif\eta\dif\xi\dif s\\
&=\frac{1}{c_1}\int^\infty_{-\infty}\int_{\mR^{2d}}\left(\int^\infty_s|\eta+\Pi^*_{s,t}\xi|^{\alpha}
\e^{-\int^t_s\psi^{\nu_r}_{\sigma_r}(-\Pi^*_{r,t}\xi-\eta)\dif r}|\hat f(t,\xi,\eta)|^2\dif t\right)\dif\eta\dif\xi\dif s\\
&=\frac{1}{c_1}\int^\infty_{-\infty}\int_{\mR^{2d}}\left(\int^t_{-\infty}|\eta+\Pi^*_{s,t}\xi|^{\alpha}
\e^{-\int^t_s\psi^{\nu_r}_{\sigma_r}(\Pi^*_{r,t}\xi+\eta)\dif r}\dif s\right)|\hat f(t,\xi,\eta)|^2\dif\eta\dif\xi\dif t\\
&\leq\left(\frac{2}{c^2_1}+\frac{2\|U\|^\alpha_\infty}{c^2_1(1+\alpha)}\right)
\int^\infty_{-\infty}\!\!\int_{\mR^{2d}}|\hat f(t,\xi,\eta)|^2\dif\eta\dif\xi\dif t
=\left(\frac{2}{c^2_1}+\frac{2\|U\|^\alpha_\infty}{c^2_1(1+\alpha)}\right)\|f\|^2_2.
\end{align*}
The proof of \eqref{BN77} for $p=2$  is thus complete by \eqref{EV10}.
\end{proof}

\end{document}